# State Estimator Design: Addressing General Delay Structures with Dissipative Constraints

Qian Feng, *Member, IEEE*, Feng Xiao, *Senior Member, IEEE*, Xiaoyu Wang, *Student Member, IEEE*

***Abstract*— Dissipative estimator (observer) design for continuous time-delay systems poses a significant challenge when an unlimited number of pointwise and general distributed delays (DDs) are concerned. We propose an effective solution to this semi-open problem using the Krasovskiĭ functional (KF) framework in conjunction with a quadratic supply rate function, where both the plant and the estimator can accommodate an unlimited number of pointwise and general distributed delays with an unlimited number of square-integrable kernels. A key contribution is the introduction of a control concept called Kronecker-Seuret Decomposition (KSD) for matrix-valued functions, which allows for the factorizations or approximations of any DD integral kernel without introducing conservatism. Moreover, using KSD facilitates the construction of complete-type KFs with integral kernels that can contain any number of weakly differentiable and linearly independent functions. Our proposed solution is formulated as sequential convex SDP problems and is set out in two theorems along with an off-line iterative algorithm, which eliminates the need for nonlinear numerical solvers. We show the effectiveness of our method using two challenging numerical experiments, including a system stabilized by a non-smooth controller.**

## I. Introduction

ESTIMATORS (Observers) [1] play a crucial role in estimating the values of state/output of dynamical systems where direct measurement of these signals is not possible. This often becomes necessary in situations wherein sensors for states or certain outputs are unavailable or unreliable, or when acquiring measurements is prohibitively costly or challenging.

The difficulties of designing estimators with performance objectives often depend on the complexity of the underlying system structures. One challenging scenario arises when systems encompass an unlimited number of pointwise and general distributed delays (DDs) in their states, where the delays can be the result of transport, propagation, and aftereffects [2]–[4] introduced by real-world media or engineering devices. This scenario can be found in various applications described by functional differential equations (FDEs) such as the modeling of wind tunnels [5], control saturation [6], event-triggered mechanisms [7], and systems with predictor controllers [8].

The primary hindrance to finding estimators for systems with general delays often stems from computational than theoretical limitations. Decades ago, Olbrot [9], [10] and Salamon [11] had obtained the sufficient and necessary condition for the existence of stable estimators for general linear time-delay systems (LTDSs) denoted in Lebesgue-Stieltjes integrals [12]. Given that the spectrum of a retarded LTDS can only contain a finite number of unstable roots [4], an estimator is always constructible by stabilizing its unstable part using methods for simple LTI systems if the system is detectable [9]–[11]. The above results, however, were quite advanced for their time, necessitating a deep understanding of spectral decomposition [13] and the computation of unstable characteristic roots, both of which present significant challenges. Other subsequent methods such as those in [14] and [15] involve solving transcendental matrix equations [16] or operator Riccati equations [15], which are also difficult to compute numerically.

The construction of Krasovskiĭ functionals (KFs) in the time-domain [17], [18] has increasingly become an effective strategy for the stability analysis and estimator design of LTDSs [19], [20], where the stability (synthesis) conditions are formulated as semidefinite programming (SDP) problems [21]. The conservatism of this approach is mainly affected by the complexity of the underlying delay systems, the generality of the predetermined forms of KFs [22], and the integral inequalities employed to construct them. A notable advantage of the KF approach is that the resulting SDP programs can be efficiently solved by various numerical algorithms [23], [24], while incorporating various control theory ingredients such as performance criteria [25] or uncertainties [26]. A comprehensive introduction to this topic can be found in [27].

To the best of our knowledge, we may state that no effective solutions to the estimator design problem of LTDSs exist in the literature when both pointwise and general distributed delays $\int_{-r}^{0} F(\tau)\boldsymbol{x}(t+\tau)\mathrm{d}\tau$ are considered along with dissipative constraints [25], [28]. Almost all existing methods on this topic only consider pointwise delays [10], [29]–[33], require special system structures [34, Eq. 12], or restrict the number of delays to one [29], [31], [33]. Since matrix kernel $F(\cdot)$ may contain general functions such that $\int_{-r}^{0} F(\tau)\mathrm{e}^{\tau s}\mathrm{d}\tau$ may not yield a closed-form expression in $s \in \mathbb{C}$, it is challenging to develop effective frameworks that can handle the complexity of general delays and dissipative constraints simultaneously in the time or frequency domains. A noteworthy approach

This work was partially supported by the National Natural Science Foundation of China under Grant Nos. 62303180, 62373150 and 62273145, Beijing Natural Science Foundation under Grant No. 4222053, and Fundamental Research Funds for Central Universities under Grant 2023MS032, China, and the ANR (France) Project Finite4SoS (ANR-15-CE23-0007). The authors are with the School of Control and Computer Engineering, North China Electric Power University Beijing, China. Email: qianfeng@ncepu.edu.cn, qfen204@aucklanduni.ac.nz, fengxiao@ncepu.edu.cn, xiaoyu_wang@ncepu.edu.cn. Corresponding Author: Feng Xiao was proposed in [20] for the $\mathcal{H}^\infty$ state estimation of an LTDS with an unlimited number of pointwise delays, using the sums-of-squares approach [35] with a coupled PDE-ODE estimator. However, when the PDE component of the estimator is implemented using various discretization techniques [36], the robustness of the discretized closed-loop system (CLS) may not always be guaranteed, which may lead to potential instability issues.

Given the preceding discourse, we propose an efficient method to design dissipative state estimators (DSEs) for a class of nonlinear TDSs using the KF approach with a quadratic supply rate function, where the final synthesis condition consists of convex SDP problems. Both the system and estimator can have an unlimited number of pointwise and general distributed delays in the states, inputs, and regulated outputs, with each DD matrix capable of having any number of $\mathcal{L}^2$ kernel functions. The structure of the system is closely aligned with the general LTDSs denoted by Lebesgue-Stieltjes integrals [12] with nonlinear components derived from output injections. Moreover, the differential equation in our study is defined in the extended sense [12] with respect to (w.r.t) the Lebesgue measure using weak derivatives, which is particularly well-suited for modeling the dynamics of engineering systems that are often exposed to noise and glitches. Owing to the generality of the formulation, our DSE design problem presents a challenge that no existing methods may handle, which could also cover a wide range of emerging engineering problems.

To address the complexity of general DDs in conjunction with the KF approach in the spirit of SDP, we propose the concept called Kronecker-Seuret Decomposition (KSD), which draws inspirations from our recent works in [19], [26] and the works by Alexandre Seuret [37] on the analysis of general DDs. Any matrix-valued functions with an unlimited number of $\mathcal{L}^2$ kernels can be decomposed by this technique as products between constant matrices and a sequence of basis functions that are linearly independent w.r.t the Lebesgue measure [38]. The vector containing these basis functions is structured with three distinct vector-valued components, where the third component is approximated by the rest of the two components using the least-squares principle [39, page 182]. A distinctive feature of the KSD concept is that it offers users the freedom to decide which kernel functions are approximated and which are directly factorized, which significantly generalizes our previous works in [19], [26]. Moreover, the expressions of KSD greatly facilitate the construction of complete type KFs whose integral kernels are the first components of the decomposition vectors. As shown by our investigations, such an arrangement allows us to increase the generality of KFs by adding an unlimited number of weakly differentiable and linearly independent functions (WDLIFs) to the integral kernels of KFs. Since the generality of KFs is primarily determined by the generality of their integral kernels, we can be assured that the synthesis conditions formulated via our KF are nonconservative. Finally, the non-conservative solvability of our framework is reinforced by the use of integral inequalities derived from the least-squares principle [39], which we employ to construct lower bounds for the quadratic integrals in our KF and its derivative.

Another hallmark of our study is the use of a functional estimator incorporating an unlimited number of pointwise delays and general DDs with $\mathcal{L}^2$ kernels, which also takes into account the effect of external disturbances under a quadratic supply rate function (SRF). The structure of our estimator is formulated in light of the general estimator utilized to describe the detectability [9], [11] of general LTDSs. Furthermore, the DDs of the estimator are constructed using the same basis functions employed in the KSD representation. Consequently, it is always possible to increase the generality of our estimator by adding new WDLIFs to the decomposition vector, potentially improving the performance of the resulting estimator.

The solution to the DSE design problem contains two theorems and an iterative algorithm, where the second theorem is derived from convexifying the bilinear matrix inequality (BMI) in the first theorem invoking the [40, Projection Lemma]. In sequence, the first theorem is computed by the proposed iterative algorithm which can be initiated by a feasible solution to the second theorem. Consequently, these three mathematical constructs act as a single package to solve the DSE synthesis problem, eliminating the need for nonlinear SDP solvers and avoiding difficult numerical schemes as seen in [14], [15].

The rest of the article is divided into four sections. In Section II, we first introduce the concept of KSD, and its utility in deriving the expression of error-dynamics using our proposed estimator. The main results concerning the dissipative estimator design are set out in Section III including two theorems and an algorithm. Next, the computation results of two numerical examples are provided in Section IV, including a challenging one involving an LTDS stabilized by a non-smooth controller. Finally, we conclude the paper and present some crucial lemmas and proofs in the appendices.

### *Notation*

$\mathbb{N}$ is the natural number set, and $\mathbb{S}^n = \{X \in \mathbb{R}^{n\times n} : X = X^\top\}$. The space of continuous functions is denoted by $\mathcal{C}(\mathcal{X};\mathbb{R}^n)$. Set $\mathcal{M}(\mathcal{X};\mathbb{R}^d)$ includes all measurable functions from $\mathcal{X}$ to $\mathbb{R}^d$. Define $\|\mathbf{x}\|^2 = \|\mathbf{x}\|_2^2 = \mathbf{x}^\top\mathbf{x}$ and $\mathcal{L}^p(\mathcal{X};\mathbb{R}^n) = \{\boldsymbol{f}(\cdot) \in \mathcal{M}(\mathcal{X};\mathbb{R}^n) : \|\boldsymbol{f}(\cdot)\|_p < +\infty\}$ with $\mathcal{X} \subseteq \mathbb{R}^n$ and the seminorm $\|\boldsymbol{f}(\cdot)\|_p^p = \int_{\mathcal{X}} \|\boldsymbol{f}(x)\|^p \mathrm{d}x$, and define Sobolev space $\mathcal{H}^1(\mathcal{X};\mathbb{R}^n) = \mathcal{W}^{1,2}(\mathcal{X};\mathbb{R}^n) = \{\boldsymbol{f}(\cdot) : \boldsymbol{f}'(\cdot) \in \mathcal{L}^2(\mathcal{X};\mathbb{R}^n)\}$, where $\boldsymbol{f}'(\cdot)$ is the weak derivative of $\boldsymbol{f}(\cdot)$. Let $\mathsf{Sy}(X) := X + X^\top$ for any $X \in \mathbb{R}^{n\times n}$. We introduce notation $\mathbf{Col}_{i=1}^n X_i = [X_i]_{i=1}^n := \left[X_1^\top \cdots X_i^\top \cdots X_n^\top\right]^\top$ to denote a column-wise concatenation of mathematical objects, whereas $\mathbf{Row}_{i=1}^n X_i = [\![X_i]\!]_{i=1}^n = \left[X_1 \cdots X_i \cdots X_n\right]$ is the row version correspondence. We also introduce notation $\widetilde{\forall} x \in \mathcal{X}$ meaning $x \in \mathcal{X}$ holds **almost everywhere** w.r.t the Lebesgue measure. Symbol $*$ is used as an abbreviations for $[*]YX = X^\top YX$ or $X^\top Y[*] = X^\top YX$ or $\left[\begin{smallmatrix} A & B \\ * & C \end{smallmatrix}\right] = \left[\begin{smallmatrix} A & B \\ B^\top & C \end{smallmatrix}\right]$. $\mathsf{O}_{n,m}$ stands for an $n\times m$ zero matrix that can be abbreviated as $\mathsf{O}_n$ if $n = m$, whereas $\mathbf{0}_n$ denotes an $n\times 1$ zero vector. We use $\oplus$ to denote $X \oplus Y = \left[\begin{smallmatrix} X & \mathsf{O}_{n,q} \\ \mathsf{O}_{p,m} & Y \end{smallmatrix}\right]$ for any $X \in \mathbb{R}^{n\times m}$, $Y \in \mathbb{R}^{p\times q}$ with its $n$-ary iterated form $\mathsf{diag}_{i=1}^\nu X_i = X_1 \oplus X_2 \oplus \cdots \oplus X_\nu$. Notation $\otimes$ is the Kronecker product. We use $\sqrt{X}$ to represent the unique square root of positive definite matrix $X \succ 0$. The order of matrix operations is defined as *matrix (scalars)* *multiplications* $>\oplus=\mathsf{diag}>\otimes>+$. Finally, we use $[\,]$ to represent empty matrices [41, I.7] follow the same definition and rules in MATLAB©. We define $I_0=[\,]_{0,0}$, $\mathsf{O}_{0,m}=[\,]_{0,m}$, $\mathsf{O}_{n,0}=[\,]_{n,0}$ and $[X_i]_{i=1}^{0}=[\,]_{0,m}$, $[\![X_i]\!]_{i=1}^{0}=[\,]_{n,0}$, where $X_i\in\mathbb{R}^{n\times m}$ and $[\,]_{0,m}$, $[\,]_{n,0}$ are empty matrices.

## II. Preliminaries and Problem Formulation

### A. Open-Loop System

Consider a TDS of the form

$$\begin{aligned}
\dot{\boldsymbol{x}}(t)=&\sum_{i=0}^{\nu}A_i\boldsymbol{x}(t-r_i)+\int_{-r}^{0}\widetilde{A}(\tau)\boldsymbol{x}(t+\tau)\mathsf{d}\tau+D_1\boldsymbol{w}(t)\\
&+\mathfrak{f}_1(\mathfrak{u}_t(\cdot),\mathfrak{y}_t(\cdot)),\qquad \forall\theta\in\mathcal{J},\ \boldsymbol{x}(t_0+\theta)=\widehat{\boldsymbol{\psi}}(\theta)\\
\boldsymbol{z}(t)=&\sum_{i=0}^{\nu}C_i\boldsymbol{x}(t-r_i)+\int_{-r}^{0}\widetilde{C}(\tau)\boldsymbol{x}(t+\tau)\mathsf{d}\tau+D_2\boldsymbol{w}(t)\\
&+\mathfrak{f}_2(\mathfrak{u}_t(\cdot),\mathfrak{y}_t(\cdot)),\qquad \forall\theta\in\mathcal{J},\ \mathfrak{u}_t(\theta)=\boldsymbol{u}(t+\theta)\\
\boldsymbol{y}(t)=&\,\mathcal{C}\boldsymbol{x}(t)+\mathfrak{f}_3(\mathfrak{u}_t(\cdot)),\quad \forall\theta\in\mathcal{J},\ \mathfrak{y}_t(\theta)=\boldsymbol{y}(t+\theta)
\end{aligned}\tag{1}$$

where the FDE holds for all $t\geq t_0\in\mathbb{R}$ almost everywhere w.r.t the Lebesgue measure. We define initial condition $\widehat{\boldsymbol{\psi}}(\cdot)\in\mathcal{C}(\mathcal{J};\mathbb{R}^n)$, and $\mathcal{J}=[-r,0]$ where the values of delays $r=r_\nu>r_{\nu-1}>\cdots>r_2>r_1>r_0=0$, $\nu\in\mathbb{N}$ are known. $\boldsymbol{x}(t)\in\mathbb{R}^n$ is the solution to the FDE in the extended sense using the Carathéodory framework [12, page 58], $\boldsymbol{u}(t)\in\mathbb{R}^p$ is the known control input, $\boldsymbol{y}(t)\in\mathbb{R}^l$ is the measurement output, $\boldsymbol{w}(\cdot)\in\mathcal{L}^2(\mathbb{R}_{\geq t_0};\mathbb{R}^q)$ denotes exogenous disturbances, and $\boldsymbol{z}(t)\in\mathbb{R}^m$ is the regulated output. Operators $\mathfrak{f}_1(\cdot)\in\mathcal{C}\left(\mathcal{C}\left(\mathcal{J};\mathbb{R}^{p+l}\right);\mathbb{R}^n\right)$, $\mathfrak{f}_2(\cdot)\in\mathcal{C}\left(\mathcal{C}\left(\mathcal{J};\mathbb{R}^{p+l}\right);\mathbb{R}^m\right)$ describe control inputs with output injection, and $\mathfrak{f}_3(\cdot)\in\mathcal{C}\left(\mathcal{C}\left(\mathcal{J};\mathbb{R}^{p}\right);\mathbb{R}^l\right)$. The dimensions of the matrices in (1) are indexed by $n;m;p;q;l\in\mathbb{N}$. Finally, the DDs in (1) satisfy

$$\widetilde{A}(\cdot)\in\mathcal{L}^2(\mathcal{J};\mathbb{R}^{n\times n}),\quad \widetilde{C}(\cdot)\in\mathcal{L}^2(\mathcal{J};\mathbb{R}^{m\times n}).$$

The integrals in (1) can always be rewritten as

$$\begin{aligned}
\int_{-r}^{0}\widetilde{A}(\tau)\boldsymbol{x}(t+\tau)\mathsf{d}\tau&=\sum_{i=1}^{\nu}\int_{\mathcal{I}_i}\widetilde{A}_i(\tau)\boldsymbol{x}(t+\tau)\mathsf{d}\tau\\
\int_{-r}^{0}\widetilde{C}(\tau)\boldsymbol{x}(t+\tau)\mathsf{d}\tau&=\sum_{i=1}^{\nu}\int_{\mathcal{I}_i}\widetilde{C}_i(\tau)\boldsymbol{x}(t+\tau)\mathsf{d}\tau
\end{aligned}\tag{2}$$

via intervals $\mathcal{I}_i=[-r_i,-r_{i-1}]$ and matrix-valued functions

$$\widetilde{A}_i(\cdot)\in\mathcal{L}^2(\mathcal{I}_i;\mathbb{R}^{n\times n}),\quad \widetilde{C}_i(\cdot)\in\mathcal{L}^2(\mathcal{I}_i;\mathbb{R}^{m\times n})\tag{3}$$

for all $i\in\mathbb{N}_\nu:=\{1,\ldots,\nu\}$. Note that decomposing $\mathcal{J}$ into $\mathcal{I}_i$ is required by the KF approach employed in the next section.

The generality of (1) is evident. Even for the linear version of (1), it covers most instances of general LTDSs described by the Lebesgue-Stieltjes integrals [12, Eq. (2), Chapter 7]. The existence of the solution of (1) is ensured by the conclusions in [12, Chapter 2], as the right-hand side of the FDE in (1) satisfies the Carathéodory conditions [12, page 58].

*Remark 1:* A wide variety of practical systems with general DDs can be modeled by (1) such as the SIR model in [42, Eq. (7)], the networked control system in [43], and the chemical reaction networks in [44, Eq. (30)]. By incorporating output injection terms $\mathfrak{f}_i(\mathfrak{u}_t(\cdot),\mathfrak{y}_t(\cdot)),i=1,2$ the nonlinear TDS in (1) can be treated as an LTDS when designing estimators. This structure can be seen as an extension of the conventional output injection scheme [45] for the case of TDSs.

It is crucial to stress that the FDE in (1) is not defined for all $t\geq t_0$ but rather for almost all $t\geq t_0$ w.r.t the Lebesgue measure [38], even with $\boldsymbol{w}(t)\equiv\boldsymbol{0}_q$. A key advantage of employing the Carathéodory formulation [12] is that the FDE in (1) is particularly well-suited for describing the dynamics of engineering and cybernetic systems, as these systems are often subject to noise and glitches with countable instances.

### B. Concept of Kronecker-Seuret Decomposition (KSD)

Since matrix-valued functions $\widetilde{A}_i(\cdot),\widetilde{C}_i(\cdot)$ in (3) have infinite dimension, directly including these terms in any SDP synthesis/stability condition that is formulated by the construction of KFs, will result in constraints of infinite dimension. To bypass this obstacle, we introduce the concept of KSD in this article, by which any DD integral matrix can be parameterized using a list of functions with matrices of finite dimensions.

*Proposition 1:* (3) holds if and only if there exist $\boldsymbol{f}_i(\cdot)\in\mathcal{H}^1(\mathcal{I}_i;\mathbb{R}^{d_i}),\boldsymbol{\varphi}_i(\cdot)\in\mathcal{L}^2(\mathcal{I}_i;\mathbb{R}^{\delta_i}),\boldsymbol{\phi}_i(\cdot)\in\mathcal{L}^2(\mathcal{I}_i;\mathbb{R}^{\mu_i})$ and $M_i\in\mathbb{R}^{d_i\times\varkappa_i}$, $\widehat{A}_i\in\mathbb{R}^{n\times\kappa_i n}$, $\widehat{C}_i\in\mathbb{R}^{m\times\kappa_i n}$ such that

$$\widetilde{A}_i(\tau)=\widehat{A}_i\left(\boldsymbol{g}_i(\tau)\otimes I_n\right),\quad \widetilde{C}_i(\tau)=\widehat{C}_i\left(\boldsymbol{g}_i(\tau)\otimes I_n\right)\tag{4}$$

$$\frac{\mathsf{d}\boldsymbol{f}_i(\tau)}{\mathsf{d}\tau}=M_i\boldsymbol{h}_i(\tau)\qquad\text{with }\ \boldsymbol{h}_i(\tau)=\begin{bmatrix}\boldsymbol{\varphi}_i(\tau)\\ \boldsymbol{f}_i(\tau)\end{bmatrix}\tag{5}$$

$$\mathsf{G}_i=\int_{\mathcal{I}_i}\boldsymbol{g}_i(\tau)\boldsymbol{g}_i^\top(\tau)\mathsf{d}\tau\succ 0\ \ \text{with }\ \boldsymbol{g}_i(\tau)=\begin{bmatrix}\boldsymbol{\phi}_i(\tau)\\ \boldsymbol{h}_i(\tau)\end{bmatrix}\tag{6}$$

hold for all $i\in\mathbb{N}_\nu$ and $\tau\in\mathcal{I}_i$, where $\kappa_i=d_i+\delta_i+\mu_i$, $\varkappa_i=d_i+\delta_i$ with $d_i\in\mathbb{N}$ and $\delta_i;\mu_i\in\mathbb{N}_0:=\mathbb{N}\cup\{0\}$, and $\frac{\mathsf{d}}{\mathsf{d}\tau}$ in (5) represents weak derivative w.r.t the Lebesgue measure.

*Proof:* First, we can infer (3) from (4)–(6) by the definitions of $\boldsymbol{\varphi}_i(\cdot),\boldsymbol{f}_i(\cdot),\boldsymbol{\phi}_i(\cdot)$ and the fact that $\mathcal{H}^1(\mathcal{I}_i;\mathbb{R}^{d_i})\subset\mathcal{L}^2(\mathcal{I}_i;\mathbb{R}^{d_i})$. This proves the necessity part of the statement.

Conversely, we seek to prove that the inclusions in (3) imply the existence of the parameters in Proposition 1 satisfying (4)–(6). Given any $\boldsymbol{f}_i(\cdot)\in\mathcal{H}^1(\mathcal{I}_i;\mathbb{R}^{d_i})$, satisfying $\forall i\in\mathbb{N}_\nu$, $\mathfrak{F}_i:=\int_{\mathcal{I}_i}\boldsymbol{f}_i(\tau)\boldsymbol{f}_i^\top(\tau)\mathsf{d}\tau\succ 0$, we can always select appropriate $\boldsymbol{\phi}_i(\cdot)\in\mathcal{L}^2(\mathcal{I}_i;\mathbb{R}^{\mu_i})$ and $\boldsymbol{\varphi}_i(\cdot)\in\mathcal{L}^2(\mathcal{I}_i;\mathbb{R}^{\delta_i})$ for some $M_i\in\mathbb{R}^{d_i\times\varkappa_i}$ to satisfy (5)–(6). This is because $\boldsymbol{f}_i'(\cdot)\in\mathcal{L}^2(\mathcal{I}_i;\mathbb{R}^{d_i})$ and $\mathsf{dim}(\boldsymbol{\varphi}_i(\tau))$ and $\mathsf{dim}(\boldsymbol{\phi}_i(\tau))$ can always be enlarged by adding new $\mathcal{L}^2$ functions that are linearly independent. Note that $\mathfrak{F}_i\succ 0$ are implied by $\mathsf{G}_i\succ 0$ in (6), indicating that the functions in $\boldsymbol{g}_i(\cdot)$ in (6) are linearly independent [46, Th. 7.2.10] in a Lebesgue sense over $\mathcal{I}_i$ for each $i\in\mathbb{N}_\nu$. We also note that also $\boldsymbol{\varphi}_i(\tau)$ or $\boldsymbol{\phi}_i(\tau)$ can be empty matrix $[\,]_{0\times 1}$.

As we can enlarge $\mathsf{dim}\left(\boldsymbol{g}_i(\tau)\right)$ in (5)–(6) without limits, we can always find appropriate $A_{i,j}\in\mathbb{R}^{n\times n}$, $C_{i,j}\in\mathbb{R}^{m\times n}$ and $\boldsymbol{g}_i^\top(\cdot)=\left[\boldsymbol{\phi}_i^\top(\cdot)\ \ \boldsymbol{\varphi}_i^\top(\cdot)\ \ \boldsymbol{f}_i^\top(\cdot)\right]=[\![g_{i,j}(\cdot)]\!]_{j=1}^{\kappa_i}$ for (3) such that

$$\widetilde{A}_i(\tau)=\sum_{j=1}^{\kappa_i}A_{i,j}g_{i,j}(\tau),\qquad \widetilde{C}_i(\tau)=\sum_{j=1}^{\kappa_i}C_{i,j}g_{i,j}(\tau)\tag{7}$$

for all $i\in\mathbb{N}_\nu$ and $\tau\in\mathcal{I}_i$ with $\kappa_i\in\mathbb{N}$, where $\boldsymbol{\varphi}_i(\cdot)\in\mathcal{L}^2(\mathcal{I}_i;\mathbb{R}^{\delta_i})$, $\boldsymbol{f}_i(\cdot)\in\mathcal{H}^1(\mathcal{I}_i;\mathbb{R}^{d_i})$ and $\boldsymbol{\phi}_i(\cdot)\in\mathcal{L}^2(\mathcal{I}_i;\mathbb{R}^{\mu_i})$

satisfy (5)–(6) for some $M_i \in \mathbb{R}^{d_i \times \varkappa_i}$. Moreover, (7) can be rewritten using the properties of matrix multiplications as

$$\widetilde{A}_i(\tau) = [\![A_{i,j}]\!]_{j=1}^{\kappa_i} G_i(\tau), \quad \widetilde{C}_i(\tau) = [\![C_{i,j}]\!]_{j=1}^{\kappa_i} G_i(\tau)$$

with $G_i(\tau) = \boldsymbol{g}_i(\tau) \otimes I_n$, which are identical to the terms in (4) by defining $\widehat{A}_i = [\![A_{i,j}]\!]_{j=1}^{\kappa_i}$ and $\widehat{C}_i = [\![C_{i,j}]\!]_{j=1}^{\kappa_i}$. ■

*Remark 2:* Since $\int_{\mathcal{I}_i} p(\tau) q(\tau) \mathsf{d}\tau$ are inner products, $\mathsf{G}_i$ in (6) are the Gramian matrices [46, Th. 7.2.10] of functions $\boldsymbol{g}_i(\cdot)$. Inequalities $\mathsf{G}_i \succ 0$ imply [46] that each row of $\boldsymbol{g}_i(\cdot)$ is linearly independent in a Lebesgue sense over $\mathcal{I}_i$ for all $i \in \mathbb{N}_\nu$. It is crucial not to confuse $\operatorname{rank} \mathsf{G}_i = \kappa_i$ with $\operatorname{rank}\left[\boldsymbol{g}_i(\tau)\boldsymbol{g}_i^\top(\tau)\right] = 1$, as $\operatorname{rank} \mathsf{G}_i \neq \operatorname{rank}\left[\boldsymbol{g}_i(\tau)\boldsymbol{g}_i^\top(\tau)\right] = 1$ due to the integrals. For a thorough understanding of the Gramian matrix for functions, see [46, Th. 7.2.10] for details. A simple example for $\mathsf{G}_i$ is the Hilbert matrices such as $\int_0^1 [\tau^i]_{i=0}^1 [\![\tau]\!]_{i=0}^1 \mathsf{d}\tau = \int_0^1 \left[\begin{smallmatrix} 1 & \tau \\ * & \tau^2 \end{smallmatrix}\right] \mathsf{d}\tau = \left[\begin{smallmatrix} 1 & 1/2 \\ * & 1/3 \end{smallmatrix}\right] \succ 0$.

Proposition 1 is the first component of the KSD concept and suffers no conservatism. A distinctive feature of the decomposition is that we can select any $\boldsymbol{f}_i(\cdot) \in \mathcal{H}^1(\mathcal{I}_i\,; \mathbb{R}^{d_i})$ to satisfy (4)–(6) even if $\widetilde{A}_i(\cdot), \widetilde{C}_i(\cdot)$ contain none of the functions in $\boldsymbol{f}_i(\cdot)$. This is because $\mathsf{dim}(\boldsymbol{f}_i(\tau)) = d_i$ and $\mathsf{dim}(\boldsymbol{\varphi}_i(\tau)) = \delta_i$ are unbounded and we can add an unlimited number of new functions to $\boldsymbol{f}_i(\cdot), \boldsymbol{\varphi}_i(\cdot)$ while still satisfying (4)–(6) with modified $\widehat{A}_i, \widehat{C}_i, M_i$. As long as all the functions in (3) are covered by some functions in $\boldsymbol{g}_i(\cdot)$, then constant matrices $\widehat{A}_i$, $\widehat{C}_i$ in (4) can always be constructed accordingly.

We separated $\boldsymbol{\varphi}_i(\cdot)$ and $\boldsymbol{\phi}_i(\cdot)$ in $\boldsymbol{g}_i(\cdot)$ as our KSD addresses them differently in its second step. Specifically, we assume $\boldsymbol{\phi}_i(\cdot)$ is approximated by $\boldsymbol{h}_i(\cdot)$ in (5) such that

$$\boldsymbol{\phi}_i(\tau) = \Gamma_i \mathfrak{H}_i^{-1} \boldsymbol{h}_i(\tau) + \boldsymbol{\varepsilon}_i(\tau), \quad \forall \tau \in \mathcal{I}_i \tag{8}$$

with $\mathfrak{H}_i = \int_{\mathcal{I}_i} \boldsymbol{h}_i(\tau) \boldsymbol{h}_i^\top(\tau) \mathsf{d}\tau$ and

$$\mathbb{R}^{\mu_i \times \varkappa_i} \ni \Gamma_i = \int_{\mathcal{I}_i} \boldsymbol{\phi}_i(\tau) \boldsymbol{h}_i^\top(\tau) \mathsf{d}\tau, \quad \forall i \in \mathbb{N}_\nu$$

using the least-square approximation [39, page 182] with $\varkappa_i, \mu_i$ in Proposition 1, where $\mathfrak{H}_i^{-1} \succ 0$ is well defined by $\mathsf{G}_i \succ 0$ in (6). Moreover, let $\boldsymbol{\varepsilon}_i(\tau) = \boldsymbol{\phi}_i(\tau) - \Gamma_i \mathfrak{H}_i^{-1} \boldsymbol{h}_i(\tau)$ denote the approximation errors and utilize matrix

$$\begin{aligned}
\mathbb{S}^{\mu_i} \ni \mathfrak{E}_i &= \int_{\mathcal{I}_i} \boldsymbol{\varepsilon}_i(\tau) \boldsymbol{\varepsilon}_i^\top(\tau) \mathsf{d}\tau \\
&= \int_{\mathcal{I}_i} \left(\boldsymbol{\phi}_i(\tau) - \Gamma_i \mathfrak{H}_i^{-1} \boldsymbol{h}_i(\tau)\right) \left(\boldsymbol{\phi}_i(\tau) - \Gamma_i \mathfrak{H}_i^{-1} \boldsymbol{h}_i(\tau)\right)^\top \mathsf{d}\tau \\
&= \int_{\mathcal{I}_i} \boldsymbol{\phi}_i(\tau) \boldsymbol{\phi}_i^\top(\tau) \mathsf{d}\tau - \Gamma_i \mathfrak{H}_i^{-1} \Gamma_i^\top, \quad \forall i \in \mathbb{N}_\nu
\end{aligned} \tag{9}$$

to quantitatively measure $\boldsymbol{\varepsilon}_i(\cdot)$, where $\mathfrak{E}_i \succ 0$ always holds on account of [19, Eq. (18)]. In view of (8), we find that

$$\begin{aligned}
\boldsymbol{g}_i(\tau) &= \begin{bmatrix} \boldsymbol{\phi}_i(\tau) \\ \boldsymbol{h}_i(\tau) \end{bmatrix} = \begin{bmatrix} \Gamma_i \mathfrak{H}_i^{-1} \boldsymbol{h}_i(\tau) \\ \boldsymbol{h}_i(\tau) \end{bmatrix} + \begin{bmatrix} \boldsymbol{\varepsilon}_i(\tau) \\ \mathbf{0}_{\varkappa_i} \end{bmatrix} \\
&= \widehat{\Gamma}_i \boldsymbol{h}_i(\tau) + \widetilde{I}_i \boldsymbol{\varepsilon}_i(\tau),
\end{aligned} \tag{10}$$

$$\widehat{\Gamma}_i = \begin{bmatrix} \Gamma_i \mathfrak{H}_i^{-1} \\ I_{\varkappa_i} \end{bmatrix} \in \mathbb{R}^{\kappa_i \times \varkappa_i}, \quad \widetilde{I}_i = \begin{bmatrix} I_{\mu_i} \\ \mathsf{O}_{\varkappa_i, \mu_i} \end{bmatrix} \in \mathbb{R}^{\kappa_i \times \mu_i}$$

with $\kappa_i, \varkappa_i, \mu_i$ in Proposition 1. With the expressions of $\boldsymbol{g}_i(\cdot)$ in (10) for (4), the KSD for $\widetilde{A}_i(\cdot), \widetilde{C}_i(\cdot)$ is established.

The structures in (4) inspire us to use "Kronecker" and "Decomposition" as a part of naming the KSD. The concept of KSD is always applicable with $\mu_i = 0$ and $\boldsymbol{\phi}_i(\cdot) = [\,]_{0\times 1}$ in Proposition 1 corresponding to the case that no functions in $\boldsymbol{g}_i(\cdot)$ are approximated, since an unlimited number of linearly independent $\mathcal{L}^2$ functions may be added to $\boldsymbol{\varphi}_i(\cdot)$. In fact, [47, Prop. 1] is a special case of the KSD concept with $\boldsymbol{\phi}_i(\cdot) = [\,]_{0\times 1}$, $\mu_i = 0$ and $\nu = 2$ without using approximations. Furthermore, the approximation scheme in [19] is also a special case of the KSD with $\delta_i = 0$ and $\nu = 1$, where the corresponding approximator $\boldsymbol{f}_1(\cdot)$ must satisfy $\frac{\mathsf{d}\boldsymbol{f}_1(\tau)}{\mathsf{d}\tau} = N\boldsymbol{f}_1(\tau)$ (standard derivatives) for some $N \in \mathbb{R}^{d \times d}$. The above condition implies that the functions in [19, $\boldsymbol{f}_1(\cdot)$] must be the solutions to linear homogeneous differential equations with constant coefficients, which are far more restrictive than of $\boldsymbol{f}_i(\cdot) \in \mathcal{H}^1(\mathcal{I}_i\,; \mathbb{R}^{d_i})$ in Proposition 1. Since the approximation scheme in [19] was inspired by the idea first pioneered by Alexandre Seuret and his colleagues in [37], where a single continuous kernel is approximated by Legendre polynomials [48, 22.2.10], we utilize "Seuret" to indicate the use of approximations in the KSD concept.

The concept of KSD represents a significant generalization of our previous work in [19], [47], both theoretically and conceptually. It amalgamates the concepts of direct decomposition [47] and least-square approximations [19] together, while substantially increasing the generality of the approximator functions $\boldsymbol{f}_i(\cdot)$. By distinguishing between $\boldsymbol{\phi}_i(\cdot)$ and $\boldsymbol{\varphi}_i(\cdot)$, it enables users to employ the approximation and factorization schemes simultaneously for the functions in $\boldsymbol{g}_i(\cdot)$ in view of their specific properties and the needs of users. This is a noted advantage of the proposed KSD, as it unifies the cases of $\boldsymbol{\phi}_i(\cdot) = [\,]_{0\times 1}$ and $\boldsymbol{\phi}_i(\cdot) \neq [\,]_{0\times 1}$ within a single framework. Additionally, the expressions and structures of KSD facilitate the utilization of general KFs in constructing synthesis conditions that encompass both the estimator gains and error matrices $\mathfrak{E}_i$. We will provide further details on this in the following sections when the main results are presented.

### *C. A Robust Estimator with General Delays*

Several decades ago, researchers [9], [49] had concluded that general LTDS $\dot{\boldsymbol{x}}(t) = \int_{-r}^0 [\mathsf{d}A(\tau)]\, \boldsymbol{x}(t+\tau) + B\boldsymbol{u}(t)$ represented by Lebesgue-Stieltjes integrals is stabilizable by state feedbacks if and only if there exists controller $\boldsymbol{u}(t) = \int_{-r}^0 [\mathsf{d}K(\tau)]\, \boldsymbol{x}(t+\tau)$ that renders the closed-loop system asymptotically stable, where the integral kernels $A(\cdot), K(\cdot)$ contain bounded variations that are right semi-continuous. As the detectability of (1) amounts to the stabilizability of its dual form [9], [11], [27], we can anticipate that estimator

$$\begin{aligned}
\dot{\widehat{\boldsymbol{x}}}(t) &= \sum_{i=0}^{\nu} A_i \widehat{\boldsymbol{x}}(t - r_i) + \int_{-r}^0 \widetilde{A}(\tau) \widehat{\boldsymbol{x}}(t+\tau) \mathsf{d}\tau \\
&\quad + \mathfrak{f}_1(\mathbf{u}_t(\cdot), \mathbf{y}_t(\cdot)) - \int_{-r}^0 [\mathsf{d}L(\tau)][\boldsymbol{y}(t+\tau) - \widehat{\boldsymbol{y}}(t+\tau)], \\
\widehat{\boldsymbol{y}}(t) &= \mathfrak{C}\widehat{\boldsymbol{x}}(t) + \mathfrak{f}_3(\mathbf{u}_t(\cdot))
\end{aligned} \tag{11}$$

achieves state estimation for (1) if and only if system (1) with $\boldsymbol{w}(t) \equiv \mathbf{0}_n$ is detectable, where $L(\tau) \in \mathbb{R}^{n \times l}$ contains bounded variations [38] that are right semi-continuous.

Inspired by the structure in (11), we consider estimator

$$\begin{aligned}
\dot{\widehat{\boldsymbol{x}}}(t) &= \mathfrak{e}(\widehat{\mathfrak{x}}_t(\cdot),\mathfrak{u}_t(\cdot),\mathfrak{y}_t(\cdot),\widehat{\mathfrak{y}}_t(\cdot))\\
&= \sum_{i=0}^{\nu} A_i\widehat{\boldsymbol{x}}(t-r_i) + \int_{-r}^{0}\widetilde{A}(\tau)\widehat{\boldsymbol{x}}(t+\tau)\mathsf{d}\tau \\
&\quad -\sum_{i=0}^{\nu} L_i[\boldsymbol{y}(t-r_i)-\widehat{\boldsymbol{y}}(t-r_i)] + \mathfrak{f}_1(\mathfrak{u}_t(\cdot),\mathfrak{y}_t(\cdot))\\
&\quad -\sum_{i=1}^{\nu}\int_{\mathcal{I}_i}\widehat{L}_i\left(\boldsymbol{g}_i(\tau)\otimes I_l\right)[\boldsymbol{y}(t+\tau)-\widehat{\boldsymbol{y}}(t+\tau)]\,\mathsf{d}\tau,\\
\widehat{\boldsymbol{z}}(t) &= \mathfrak{z}(\widehat{\mathfrak{x}}_t(\cdot),\mathfrak{u}_t(\cdot),\mathfrak{y}_t(\cdot),\widehat{\mathfrak{y}}_t(\cdot)) \qquad\qquad (12)\\
&= \sum_{i=0}^{\nu} C_i\widehat{\boldsymbol{x}}(t-r_i) + \int_{-r}^{0}\widetilde{C}(\tau)\widehat{\boldsymbol{x}}(t+\tau)\mathsf{d}\tau \\
&\quad -\sum_{i=0}^{\nu} \mathcal{L}_i[\boldsymbol{y}(t-r_i)-\widehat{\boldsymbol{y}}(t-r_i)] + \mathfrak{f}_2(\mathfrak{u}_t(\cdot),\mathfrak{y}_t(\cdot))\\
&\quad -\sum_{i=1}^{\nu}\int_{\mathcal{I}_i}\widehat{\mathcal{L}}_i\left(\boldsymbol{g}_i(\tau)\otimes I_l\right)[\boldsymbol{y}(t+\tau)-\widehat{\boldsymbol{y}}(t+\tau)]\,\mathsf{d}\tau,\\
\widehat{\boldsymbol{y}}(t) &= \mathcal{C}\widehat{\boldsymbol{x}}(t) + \mathfrak{f}_3(\mathfrak{u}_t(\cdot)),\\
&\forall\theta\in\mathcal{J},\ \widehat{\mathfrak{x}}_t(\theta)=\widehat{\boldsymbol{x}}(t+\theta),\quad \widehat{\mathfrak{y}}_t(\theta)=\widehat{\boldsymbol{y}}(t+\theta)
\end{aligned}$$

for estimating $\boldsymbol{x}(t)$ in (1) with $\boldsymbol{g}_i(\cdot)$ in Proposition 1, where $L_i \in \mathbb{R}^{n\times l}$, $\widehat{L}_i \in \mathbb{R}^{n\times l\kappa_i}$ and $\mathcal{L}_i \in \mathbb{R}^{m\times l}$, $\widehat{\mathcal{L}}_i \in \mathbb{R}^{m\times l\kappa_i}$ are the gains yet to be determined, and $\widehat{\boldsymbol{x}}(t), \widehat{\boldsymbol{z}}(t), \widehat{\boldsymbol{y}}(t)$ are the estimated signals for $\boldsymbol{x}(t), \boldsymbol{z}(t), \boldsymbol{y}(t)$, respectively. The structure in (12) can be seen as a special case of (11), tailored to accommodate the complex delay structures in (1) using the KSD concept. In fact, it can be regarded as an example of the Luenberger-type of observers for TDSs. A notable feature of (12) is that $\widehat{\boldsymbol{z}}(t)$ is also compensated by measurement output injection, independent of the estimator gains for $\boldsymbol{x}(t)$. This design could improve the performance of estimating $\boldsymbol{z}(t)$ as a by-product with respect to various objectives.

*Remark 3:* Though the use of (11) avoids conservatism, it is nearly impossible to compute the Lebesgue-Stieltjes integral kernel $L(\tau)$ numerically. Consequently, we propose estimator (12) to ensure that the gains can be computed via an application of the KSD concept without significantly compromising the generality of the estimator compared to (11).

Our goal is to find appropriate $L_i$, $\widehat{L}_i$, $\mathcal{L}_i$, $\widehat{\mathcal{L}}_i$ to achieve $\lim_{t\to\infty}\|\boldsymbol{x}(t)-\widehat{\boldsymbol{x}}(t)\| = 0$ with $\boldsymbol{w}(t) \equiv \mathbf{0}_q$, which also ensures $\lim_{t\to\infty}\|\boldsymbol{z}(t)-\widehat{\boldsymbol{z}}(t)\| = 0$ with $\boldsymbol{w}(t)\equiv\mathbf{0}_q$ as a by-product. Meanwhile, it is more reasonable to assume that the operation of estimator (12) is affected by external disturbances from the operating environment. Thus, we use

$$\begin{aligned}
\dot{\widehat{\boldsymbol{x}}}(t) &= \mathfrak{e}(\widehat{\mathfrak{x}}_t(\cdot),\mathfrak{u}_t(\cdot),\mathfrak{y}_t(\cdot),\widehat{\mathfrak{y}}_t(\cdot)) + D_3\boldsymbol{w}(t),\\
\widehat{\boldsymbol{z}}(t) &= \mathfrak{z}(\widehat{\mathfrak{x}}_t(\cdot),\mathfrak{u}_t(\cdot),\mathfrak{y}_t(\cdot),\widehat{\mathfrak{y}}_t(\cdot)) + D_4\boldsymbol{w}(t), \qquad (13)\\
\widehat{\boldsymbol{y}}(t) &= \mathcal{C}\widehat{\boldsymbol{x}}(t) + \mathfrak{f}_3(\mathfrak{u}_t(\cdot))
\end{aligned}$$

to denote the dynamics in such contexts, where $D_3\boldsymbol{w}(t)$ and $D_4\boldsymbol{w}(t)$ represent disturbances from the operating environment or modeling errors. It is crucial to note that $D_3\boldsymbol{w}(t)$ and $D_4\boldsymbol{w}(t)$ are not part of estimator (12) as they originate from the operating environment. Rather, we utilize (13) in the subsequent analyses to ensure that (12) has strong robustness. Once the gains in (12) are obtained, estimator (12) is implemented without requiring the measurement of $\boldsymbol{w}(t)$.

## D. Formulation of Error-Dynamics

We give the formulation of the error-dynamics in this subsection using (1) and (13) via an application of the KSD concept, where the following lemma is extensively utilized in the derivations. We kindly advise readers to be familiar with the matrix identities in this lemma before proceeding.

*Lemma 1:* $\forall X\in\mathbb{R}^{n\times m},\ \forall Y\in\mathbb{R}^{m\times p},\ \forall Z\in\mathbb{R}^{q\times r}$,

$$\begin{aligned}
(X\otimes I_q)(Y\otimes Z) &= XY\otimes Z\\
&= XY\otimes ZI_r = (X\otimes Z)(Y\otimes I_r)\\
&= I_m XY\otimes(ZI_r) = (I_m\otimes Z)(XY\otimes I_r). \qquad (14a)
\end{aligned}$$

Moreover, $\forall X\in\mathbb{R}^{n\times m}$, we see that

$$\begin{bmatrix} A & B\\ C & D\end{bmatrix}\otimes X = \begin{bmatrix} A\otimes X & B\otimes X\\ C\otimes X & D\otimes X\end{bmatrix},\ I_n\otimes X = \operatorname*{diag}_{i=1}^{n} X \qquad (14b)$$

for any $A, B, C, D$ with appropriate dimensions. Finally, for any matrix $X_i, Y_i, Z_i$ and vector $\mathbf{v}_i$ with matching dimensions,

$$\begin{aligned}
\sum_{i=1}^{n}(X_i+Y_iZ_i)\mathbf{v}_i &= [\![X_i+Y_iZ_i]\!]_{i=1}^{n}\,[\mathbf{v}_i]_{i=1}^{n}\\
&= \left([\![X_i]\!]_{i=1}^{n} + [\![Y_iZ_i]\!]_{i=1}^{n}\right)[\mathbf{v}_i]_{i=1}^{n}\\
&= [\![X_i]\!]_{i=1}^{n}\,[\mathbf{v}_i]_{i=1}^{n} + [\![Y_iZ_i]\!]_{i=1}^{n}\,[\mathbf{v}_i]_{i=1}^{n}\\
&= [\![X_i]\!]_{i=1}^{n}\,[\mathbf{v}_i]_{i=1}^{n} + [\![Y_i]\!]_{i=1}^{n}\left[\operatorname*{diag}_{i=1}^{n} Z_i\right][\mathbf{v}_i]_{i=1}^{n}, \qquad (14c)
\end{aligned}$$

$$\begin{aligned}
\sum_{i=1}^{n} Y_iZ_i\mathbf{v}_i &= [\![Y_i]\!]_{i=1}^{n}\left[\operatorname*{diag}_{i=1}^{n} Z_i\right][\mathbf{v}_i]_{i=1}^{n}\\
&= \begin{bmatrix}[\![Y_i]\!]_{i=1}^{n} & \mathsf{O}\end{bmatrix}\left[\left[\operatorname*{diag}_{i=1}^{n} Z_i\right]\oplus\mathsf{O}\right][\mathbf{v}_i]_{i=1}^{n+1}, \qquad (14d)
\end{aligned}$$

$$[Y_iZ_i]_{i=1}^{n} = \left(\operatorname*{diag}_{i=1}^{n} Y_i\right)[Z_i]_{i=1}^{n}, \qquad (14e)$$

$$\begin{aligned}
\operatorname*{diag}_{i=1}^{n} Y_iZ_i\otimes I_n &= \left(\left[\operatorname*{diag}_{i=1}^{n} Y_i\right]\left[\operatorname*{diag}_{i=1}^{n} Z_i\right]\right)\otimes I_n\\
&= \left[\operatorname*{diag}_{i=1}^{n} Y_i\otimes I_n\right]\left[\operatorname*{diag}_{i=1}^{n} Z_i\otimes I_n\right] = \operatorname*{diag}_{i=1}^{n}\left(Y_iZ_i\otimes I_n\right). \quad (14f)
\end{aligned}$$

*Proof:* The proof is obtained by the properties of matrix multiplications and the Kronecker product [46] together with the definition of $n$-ary operators $\mathrm{diag}_{i=1}^{n}$, $[\![\cdot]\!]_{i=1}^{n}$ and $[\cdot]_{i=1}^{n}$. ■

Let $\boldsymbol{e}(t) = \boldsymbol{x}(t)-\widehat{\boldsymbol{x}}(t)$ and $\boldsymbol{\zeta}(t) = \boldsymbol{z}(t)-\widehat{\boldsymbol{z}}(t)$ be the estimation error for $\boldsymbol{x}(t)$ and $\boldsymbol{z}(t)$, respectively, and

$$\boldsymbol{\chi}(t,\theta) = [\boldsymbol{e}(t+\acute{r}_i\theta-r_{i-1})]_{i=1}^{\nu}\in\mathbb{R}^{\nu n},\ \theta\in[-1,0] \qquad (15)$$

with $\acute{r}_i = r_i - r_{i-1}$, whose structure is inspired by the state variable of the PDE in [50]. By using (14a) in light of $\boldsymbol{y}(t)$ in (1) and $\widehat{\boldsymbol{y}}(t)$ in (12) with $\boldsymbol{g}_i(\tau)$ in Proposition 1, then

$$\begin{aligned}
\forall i\in\mathbb{N}_\nu,\ &(\boldsymbol{g}_i(\tau)\otimes I_l)\,[\boldsymbol{y}(t+\tau)-\widehat{\boldsymbol{y}}(t+\tau)]\\
&= (\boldsymbol{g}_i(\tau)\otimes I_l)\,\mathcal{C}\boldsymbol{e}(t+\tau) = (I_{\kappa_i}\otimes\mathcal{C})G_i(\tau)\boldsymbol{e}(t+\tau)
\end{aligned}$$

with $G_i(\tau) = \boldsymbol{g}_i(\tau)\otimes I_n$ and $\kappa_i$ in Proposition 1, and

$$\sum_{i=1}^{\nu}\int_{\mathcal{I}_i}\widehat{L}_i\left(\boldsymbol{g}_i(\tau)\otimes I_l\right)[\boldsymbol{y}(t+\tau)-\widehat{\boldsymbol{y}}(t+\tau)]\,\mathsf{d}\tau$$

$$
= \sum_{i=1}^{\nu} \int_{\mathcal{I}_i} \widehat{L}_i (I_{\kappa_i} \otimes \mathfrak{C}) G_i(\tau) \boldsymbol{e}(t+\tau) \mathsf{d}\tau, \tag{16a}
$$

$$
\sum_{i=1}^{\nu} \int_{\mathcal{I}_i} \widehat{\mathfrak{L}}_i \left(\boldsymbol{g}_i(\tau) \otimes I_l\right) \left[\boldsymbol{y}(t+\tau) - \widehat{\boldsymbol{y}}(t+\tau)\right] \mathsf{d}\tau
$$

$$
= \sum_{i=1}^{\nu} \int_{\mathcal{I}_i} \widehat{\mathfrak{L}}_i (I_{\kappa_i} \otimes \mathfrak{C}) G_i(\tau) \boldsymbol{e}(t+\tau) \mathsf{d}\tau. \tag{16b}
$$

Now employ Proposition 1 to the DDs in (1) with (2)–(3) and (16a)–(16b), and then substrate (13), this leads to

$$
\begin{aligned}
\dot{\boldsymbol{e}}(t) &= \llbracket A_i + L_i \mathfrak{C} \rrbracket_{i=0}^{\nu} \begin{bmatrix} \boldsymbol{e}^\top(t) & \boldsymbol{\chi}^\top(t,-1) \end{bmatrix}^\top \\
&\quad + \sum_{i=1}^{\nu} \int_{\mathcal{I}_i} \left(\widehat{A}_i + \widehat{L}_i (I_{\kappa_i} \otimes \mathfrak{C})\right) G_i(\tau) \boldsymbol{e}(t+\tau) \mathsf{d}\tau \\
&\quad + (D_1 - D_3)\, \boldsymbol{w}(t), \\
\boldsymbol{\zeta}(t) &= \llbracket (C_i + \mathfrak{L}_i \mathfrak{C}) \rrbracket_{i=0}^{\nu} \begin{bmatrix} \boldsymbol{e}^\top(t) & \boldsymbol{\chi}^\top(t,-1) \end{bmatrix}^\top \\
&\quad + \sum_{i=1}^{\nu} \int_{\mathcal{I}_i} \left(\widehat{C}_i + \widehat{\mathfrak{L}}_i (I_{\kappa_i} \otimes \mathfrak{C})\right) G_i(\tau) \boldsymbol{e}(t+\tau) \mathsf{d}\tau \\
&\quad + (D_2 - D_4)\, \boldsymbol{w}(t), \\
&\forall \theta \in \mathcal{J}, \ \ \boldsymbol{e}(t_0 + \theta) = \boldsymbol{\psi}(\tau), \ \ \boldsymbol{\psi}(\cdot) \in \mathcal{C}(\mathcal{J}; \mathbb{R}^n)
\end{aligned} \tag{17}
$$

by the property in (14c) with $\boldsymbol{\chi}(t,-1) = [\boldsymbol{e}(t - r_i)]_{i=1}^{\nu}$ .

*Remark 4:* The error-dynamics in (17) is independent of the control input/output feedback $\mathfrak{f}_i(\mathfrak{u}_t(\cdot), \mathfrak{y}_t(\cdot))$ in (1). This demonstrates that the state/output estimation problem of nonlinear TDS (1) can be addressed similar to LTDSs by utilizing the estimator structure in (12).

*Remark 5:* Since the error-dynamics in (17) is a retarded LTDS that has no neutral terms [4], its stability with $\boldsymbol{w}(t) \equiv \mathbf{0}_q$ is always preserved when the DDs in (12) are numerically approximated using ordinary quadratures such as the trapezoidal rule, as long as the accuracy of approximations reaches certain levels [51, Th. 2]. This ensures the safe implementation of (12) in practical settings without requiring special approximation rules [52] for the DD integrals.

Now by the expressions in (10) with (14a), we find that

$$
\begin{aligned}
G_i(\tau) = \boldsymbol{g}_i(\tau) \otimes I_n &= \left[\widehat{\Gamma}_i \boldsymbol{h}_i(\tau) + \widetilde{I}_i \boldsymbol{\varepsilon}_i(\tau)\right] \otimes I_n \\
&= \left(\widehat{\Gamma}_i \otimes I_n\right) H_i(\tau) + \left(\widetilde{I}_i \otimes I_n\right) E_i(\tau)
\end{aligned} \tag{18}
$$

$$
\text{and } (I_{\kappa_i} \otimes \mathfrak{C}) G_i(\tau) = \left(\widehat{\Gamma}_i \otimes \mathfrak{C}\right) H_i(\tau) + \left(\widetilde{I}_i \otimes \mathfrak{C}\right) E_i(\tau),
$$

$$
\text{where } H_i(\tau) = \boldsymbol{h}_i(\tau) \otimes I_n \ \text{ and } \ E_i(\tau) = \boldsymbol{\varepsilon}_i(\tau) \otimes I_n. \tag{19}
$$

By utilizing (14a) with (18)–(19) and $\mathfrak{E}_i \succ 0$, $\mathfrak{H}_i \succ 0$ in (8)–(9), we rewrite the DD integral matrices in (17) as

$$
\begin{aligned}
&\left[\widehat{A}_i + \widehat{L}_i (I_{\kappa_i} \otimes \mathfrak{C})\right] G_i(\tau) \\
&= \left[\widehat{A}_i \left(\widehat{\Gamma}_i \otimes I_n\right) + \widehat{L}_i \left(\widehat{\Gamma}_i \otimes \mathfrak{C}\right)\right] H_i(\tau) \\
&\qquad + \left[\widehat{A}_i \left(\widetilde{I}_i \otimes I_n\right) + \widehat{L}_i \left(\widetilde{I}_i \otimes \mathfrak{C}\right)\right] E_i(\tau) \\
&= \left[\widehat{A}_i (T_i \otimes I_n) + \widehat{L}_i (T_i \otimes \mathfrak{C})\right] \left[\sqrt{\mathfrak{H}_i^{-1}} \otimes I_n\right] H_i(\tau) \\
&\quad + \left[\widehat{A}_i \left(\widetilde{T}_i \otimes I_n\right) + \widehat{L}_i \left(\widetilde{T}_i \otimes \mathfrak{C}\right)\right] \left[\sqrt{\mathfrak{E}_i^{-1}} \otimes I_n\right] E_i(\tau)
\end{aligned} \tag{20}
$$

$$
\begin{aligned}
&\text{and } \left[\widehat{C}_i + \widehat{\mathfrak{L}}_i (I_{\kappa_i} \otimes \mathfrak{C})\right] G_i(\tau) \\
&= \left[\widehat{C}_i (T_i \otimes I_n) + \widehat{\mathfrak{L}}_i (T_i \otimes \mathfrak{C})\right]\left[\sqrt{\mathfrak{H}_i^{-1}} \otimes I_n\right] H_i(\tau) \\
&\quad + \left[\widehat{C}_i \left(\widetilde{T}_i \otimes I_n\right) + \widehat{\mathfrak{L}}_i \left(\widetilde{T}_i \otimes \mathfrak{C}\right)\right]\left[\sqrt{\mathfrak{E}_i^{-1}} \otimes I_n\right] E_i(\tau)
\end{aligned} \tag{21}
$$

for all $i \in \mathbb{N}_\nu$, where $T_i = \widehat{\Gamma}_i \sqrt{\mathfrak{H}_i} = \begin{bmatrix} \Gamma_i \sqrt{\mathfrak{H}_i^{-1}} \\ \sqrt{\mathfrak{H}_i} \end{bmatrix} \in \mathbb{R}^{\kappa_i \times \varkappa_i}$ and $\widetilde{T}_i = \widetilde{I}_i \sqrt{\mathfrak{E}_i}$. We kindly remind readers that the order of matrix operations was defined in the introduction.

By making use of the crucial properties in Lemma 1 again with (20)–(21), we can reformulate the terms in (17) as

$$
\begin{aligned}
\llbracket A_i + L_i \mathfrak{C} \rrbracket_{i=0}^{\nu} &= \llbracket A_i \rrbracket_{i=0}^{\nu} + \llbracket L_i \rrbracket_{i=0}^{\nu} \left(I_{\nu+1} \otimes \mathfrak{C}\right), \\
\llbracket C_i + \mathfrak{L}_i \mathfrak{C} \rrbracket_{i=0}^{\nu} &= \llbracket C_i \rrbracket_{i=0}^{\nu} + \llbracket \mathfrak{L}_i \rrbracket_{i=0}^{\nu} \left(I_{\nu+1} \otimes \mathfrak{C}\right),
\end{aligned} \tag{22}
$$

$$
\begin{aligned}
&\text{and } \sum_{i=1}^{\nu} \int_{\mathcal{I}_i} \left[\widehat{A}_i + \widehat{L}_i \left(I_{\kappa_i} \otimes \mathfrak{C}\right)\right] G_i(\tau) \boldsymbol{e}(t+\tau) \mathsf{d}\tau \\
&= \left\llbracket \widehat{A}_i (T_i \otimes I_n) \right\rrbracket_{i=1}^{\nu} \boldsymbol{\xi}(t) + \left\llbracket \widehat{L}_i (T_i \otimes I_n) \right\rrbracket_{i=1}^{\nu} (I_{\varkappa\nu} \otimes \mathfrak{C}) \boldsymbol{\xi}(t) \\
&+ \left\llbracket \widehat{A}_i \left(\widetilde{T}_i \otimes I_n\right) \right\rrbracket_{i=1}^{\nu} \mathfrak{a}(t) + \left\llbracket \widehat{L}_i \left(\widetilde{T}_i \otimes I_n\right) \right\rrbracket_{i=1}^{\nu} (I_{\mu\nu} \otimes \mathfrak{C}) \mathfrak{a}(t),
\end{aligned}
$$

$$
\begin{aligned}
&\text{and } \sum_{i=1}^{\nu} \int_{\mathcal{I}_i} \left[\widehat{C}_i + \widehat{\mathfrak{L}}_i \left(I_{\kappa_i} \otimes \mathfrak{C}\right)\right] G_i(\tau) \boldsymbol{e}(t+\tau) \mathsf{d}\tau \\
&= \left\llbracket \widehat{C}_i (T_i \otimes I_n) \right\rrbracket_{i=1}^{\nu} \boldsymbol{\xi}(t) + \left\llbracket \widehat{\mathfrak{L}}_i (T_i \otimes I_n) \right\rrbracket_{i=1}^{\nu} (I_{\varkappa\nu} \otimes \mathfrak{C}) \boldsymbol{\xi}(t) \\
&+ \left\llbracket \widehat{C}_i \left(\widetilde{T}_i \otimes I_n\right) \right\rrbracket_{i=1}^{\nu} \mathfrak{a}(t) + \left\llbracket \widehat{\mathfrak{L}}_i \left(\widetilde{T}_i \otimes I_n\right) \right\rrbracket_{i=1}^{\nu} (I_{\mu\nu} \otimes \mathfrak{C}) \mathfrak{a}(t),
\end{aligned} \tag{23}
$$

where $\varkappa = \sum_{i=1}^{\nu} \varkappa_i$ and $\mu = \sum_{i=1}^{\nu} \mu_i$ and

$$
\begin{aligned}
\boldsymbol{\xi}(t) &= \left[\int_{\mathcal{I}_i} \left(\sqrt{\mathfrak{H}_i^{-1}} \otimes I_n\right) H_i(\tau) \boldsymbol{e}(t+\tau) \mathsf{d}\tau\right]_{i=1}^{\nu}, \\
\mathfrak{a}(t) &= \left[\int_{\mathcal{I}_i} \left(\sqrt{\mathfrak{E}_i^{-1}} \otimes I_n\right) E_i(\tau) \boldsymbol{e}(t+\tau) \mathsf{d}\tau\right]_{i=1}^{\nu}
\end{aligned} \tag{24}
$$

with $H_i(\tau), E_i(\tau)$ in (19) and $\varkappa_i, \mu_i$ in Proposition 1. Finally, applying (22)–(24) to (17) with Lemma 1 produces

$$
\begin{aligned}
\dot{\boldsymbol{e}}(t) &= \left(\mathbf{A} + \mathbf{L}_1 \left[(I_\beta \otimes \mathfrak{C}) \oplus \mathsf{O}_q\right]\right) \boldsymbol{\vartheta}(t), \ \ \widetilde{\forall} t \geq t_0 \in \mathbb{R} \\
\boldsymbol{\zeta}(t) &= \left(\mathbf{C} + \mathbf{L}_2 \left[(I_\beta \otimes \mathfrak{C}) \oplus \mathsf{O}_q\right]\right) \boldsymbol{\vartheta}(t), \\
&\forall \theta \in \mathcal{J}, \ \boldsymbol{e}(t_0 + \theta) = \boldsymbol{\psi}(\theta), \ \boldsymbol{\psi}(\cdot) \in \mathcal{C}\left(\mathcal{J}; \mathbb{R}^n\right)
\end{aligned} \tag{25}
$$

with $\beta = 1 + \nu + \kappa$ and $\kappa = \varkappa + \mu = \sum_{i=1}^{\nu} \kappa_i$, where

$$
\boldsymbol{\vartheta}(t) = \begin{bmatrix} \boldsymbol{\omega}^\top(t) & \mathfrak{a}^\top(t) & \boldsymbol{w}^\top(t) \end{bmatrix}^\top, \tag{26a}
$$

$$
\boldsymbol{\omega}(t) = \begin{bmatrix} \boldsymbol{e}^\top(t) & \boldsymbol{\chi}^\top(t,-1) & \boldsymbol{\xi}^\top(t) \end{bmatrix}^\top, \tag{26b}
$$

$$
\begin{aligned}
\mathbf{A} = \Big[ A_0 \quad \llbracket A_i \rrbracket_{i=1}^{\nu} \quad \left\llbracket \widehat{A}_i (T_i \otimes I_n) \right\rrbracket_{i=1}^{\nu} \cdots \\
\cdots \left\llbracket \widehat{A}_i \left(\widetilde{T}_i \otimes I_n\right) \right\rrbracket_{i=1}^{\nu} \quad D_1 - D_3 \Big],
\end{aligned} \tag{26c}
$$

$$
\begin{aligned}
\mathbf{C} = \Big[ C_0 \quad \llbracket C_i \rrbracket_{i=1}^{\nu} \quad \left\llbracket \widehat{C}_i (T_i \otimes I_n) \right\rrbracket_{i=1}^{\nu} \cdots \\
\cdots \left\llbracket \widehat{C}_i \left(\widetilde{T}_i \otimes I_n\right) \right\rrbracket_{i=1}^{\nu} \quad D_2 - D_4 \Big],
\end{aligned} \tag{26d}
$$

$$
\begin{aligned}
\mathbf{L}_1 = \Big[ L_0 \quad \llbracket L_i \rrbracket_{i=1}^{\nu} \quad \left\llbracket \widehat{L}_i (T_i \otimes I_l) \right\rrbracket_{i=1}^{\nu} \cdots \\
\cdots \left\llbracket \widehat{L}_i \left(\widetilde{T}_i \otimes I_l\right) \right\rrbracket_{i=1}^{\nu} \quad \mathsf{O}_{n,q} \Big],
\end{aligned} \tag{26e}
$$

$$\mathbf{L}_2 = \Big[\mathcal{L}_0 \quad [\![\mathcal{L}_i]\!]_{i=1}^{\nu} \quad \Big[\!\Big[\widehat{\mathcal{L}}_i \left(T_i \otimes I_l\right)\Big]\!\Big]_{i=1}^{\nu} \cdots$$
$$\cdots \Big[\!\Big[\widehat{\mathcal{L}}_i \left(\widetilde{T}_i \otimes I_l\right)\Big]\!\Big]_{i=1}^{\nu} \qquad \mathsf{O}_{m,q}\Big]. \quad (26f)$$

Similar to the formulation in (1), the FDE in (25) holds for almost all $t \geq t_0$ w.r.t the Lebesgue measure, where the existence and uniqueness [12, page 58] of its solution are guaranteed, given that the system is linear in $\boldsymbol{e}(\cdot)$.

## III. Main Results on the DSE design problem

We first give some propositions that can characterize the stability and dissipativity of (25). Since the FDE in (25) is defined in the extended sense which means that the right hand side of (25) may not be piecewise continuous in $t$, the Krasovskiĭ stability theorems [12], [53] for FDEs with classical derivatives may not be applied to the stability of (25). Thus, we have developed a variant of Krasovskiĭ stability criterion for the analysis of FDEs in the extended sense.

*Lemma 2:* Let $\boldsymbol{w}(t) \equiv \mathbf{0}_q$ in (25) and all delay values $r_i > 0$ be given. Then the trivial solution $\boldsymbol{e}(t) \equiv \mathbf{0}_n$ of (25) is exponentially stable for any $\boldsymbol{\psi}(\cdot) \in \mathcal{C}(\mathcal{J};\mathbb{R}^n)$ in (25) if there exist $\epsilon_1;\epsilon_2;\epsilon_3 > 0$ and $\mathsf{v}(\cdot) \in \mathcal{C}^1(\mathcal{C}(\mathcal{J};\mathbb{R}^n);\mathbb{R})$ such that

$$\epsilon_1 \left\|\boldsymbol{\psi}(0)\right\|^2 \leq \mathsf{v}(\boldsymbol{\psi}(\cdot)) \leq \epsilon_2 \left\|\boldsymbol{\psi}(\cdot)\right\|_{\infty}^2, \qquad (27a)$$

$$\widetilde{\forall} t \geq t_0, \ \tfrac{\mathsf{d}}{\mathsf{d}t}\mathsf{v}(\boldsymbol{e}_t(\cdot)) \leq -\,\epsilon_3 \left\|\boldsymbol{e}(t)\right\|^2 \qquad (27b)$$

for any function $\boldsymbol{\psi}(\cdot) \in \mathcal{C}(\mathcal{J};\mathbb{R}^n)$ in (25), where $\left\|\boldsymbol{\psi}(\cdot)\right\|_{\infty}^2 := \sup_{\tau\in\mathcal{J}} \left\|\boldsymbol{\psi}(\tau)\right\|^2$. Furthermore, $\boldsymbol{e}_t(\cdot)$ in (27b) is defined by $\forall t \geq t_0$, $\forall \theta \in \mathcal{J}$, $\boldsymbol{e}_t(\theta) = \boldsymbol{e}(t+\theta)$ in which $\boldsymbol{e} : [t_0 - r, \infty) \to \mathbb{R}^n$ satisfies (25) with $\boldsymbol{w}(t) \equiv \mathbf{0}_q$.

*Proof:* The proof is similar to that of [47, Cor. 1] which is a special case of [47, Lem. 4]. Since (25) is linear and retarded, asymptotic stability ensures exponential stability. ■

We require the following characterization of dissipativity.

*Definition 1:* System (25) with a supply rate function (SRF) $\mathsf{s}(\boldsymbol{\zeta}(t),\boldsymbol{w}(t))$ is dissipative if there exists a differentiable functional $\mathsf{v} : \mathcal{C}(\mathcal{J};\mathbb{R}^n) \to \mathbb{R}$ such that

$$\widetilde{\forall} t \geq t_0, \quad \tfrac{\mathsf{d}}{\mathsf{d}t}\mathsf{v}(\boldsymbol{e}_t(\cdot)) - \mathsf{s}(\boldsymbol{\zeta}(t),\boldsymbol{w}(t)) \leq 0 \qquad (28)$$

with $t_0 \in \mathbb{R}$ and $\boldsymbol{\zeta}(t)$ in (25). Moreover, $\boldsymbol{e}_t(\cdot)$ in (28) is defined by $\forall t \geq t_0, \forall \theta \in \mathcal{J}, \boldsymbol{e}_t(\theta) = \boldsymbol{e}(t+\theta)$ with $\boldsymbol{e}(t)$ satisfying (25). If (28) holds, then it follows that

$$\forall t \geq t_0, \ \ \mathsf{v}(\boldsymbol{e}_t(\cdot)) - \mathsf{v}(\boldsymbol{e}_{t_0}(\cdot)) \leq \int_{t_0}^{t} \mathsf{s}(\boldsymbol{z}(\theta),\boldsymbol{w}(\theta))\mathsf{d}\theta \quad (29)$$

from the fundamental theorem [38] of Lebesgue integrals, which is the original definition of dissipativity [28]. This is because $\dot{\mathsf{v}}(\boldsymbol{x}_t(\cdot))$ is integrable and exists for almost all $t \geq t_0$.

In this paper, we utilize a quadratic SRF

$$\mathsf{s}(\boldsymbol{\zeta}(t),\boldsymbol{w}(t)) = \begin{bmatrix}\boldsymbol{\zeta}(t)\\ \boldsymbol{w}(t)\end{bmatrix}^{\top}\begin{bmatrix}\widetilde{J}^{\top}J_1^{-1}\widetilde{J} & J_2\\ * & J_3\end{bmatrix}\begin{bmatrix}\boldsymbol{\zeta}(t)\\ \boldsymbol{w}(t)\end{bmatrix}, \qquad (30)$$

$$\widetilde{J}^{\top}J_1^{-1}\widetilde{J} \preceq 0, \ J_1^{-1} \prec 0, \ \widetilde{J} \in \mathbb{R}^{m\times m}, \ J_2 \in \mathbb{R}^{m\times q}, \ J_3 \in \mathbb{S}^q$$

to enforce dissipativity, where matrices $\widetilde{J}, J_2$ are known. The structure in (30) is capable of describing numerous performance criteria [25] such as

- $\mathcal{L}^2$ gain performance: $J_1 = -\gamma I_m$, $\widetilde{J} = I_m$, $J_2 = \mathsf{O}_{m,q}$, $J_3 = \gamma I_q$ with free variable $\gamma > 0$;
- Strict Passivity: $J_1 \prec 0$, $\widetilde{J} = \mathsf{O}_m$, $J_2 = I_m$, $J_3 = \mathsf{O}_m$;
- Other sector constraints in [54, Table 1].

Aiming to facilitate our mathematical derivations, the SRF in (30) has minor differences from the one in [25] which remains mathematically equivalent.

With error-dynamics (25), the main results in this paper are presented in Theorems 1–2 and Algorithm 1, where the proofs and remarks of Theorems 1–2 are **located in the Appendix.** Moreover, a dedicated analysis of the conservatism of Theorem 1 is provided in Appendix B right after its proof.

*Theorem 1:* Let $\boldsymbol{f}_i(\cdot), \boldsymbol{\varphi}_i(\cdot), \boldsymbol{\phi}_i(\cdot), M_i$ in Proposition 1 and $\Gamma_i, \mathfrak{F}_i, \mathfrak{H}_i, \mathfrak{E}_i$ in (8)–(10) be given for an application of the KSD concept in Section II-B. Then the error-dynamics in (25) with SRF (30) is dissipative, and the trivial solution to (25) with $\boldsymbol{w}(t) \equiv \mathbf{0}_q$ is exponentially stable for any $\boldsymbol{\psi}(\cdot) \in \mathcal{C}(\mathcal{J};\mathbb{R}^n)$ in (25) if there exist estimator gains $L_0; L_i \in \mathbb{R}^{n\times l}$, $\widehat{L}_i \in \mathbb{R}^{n\times\kappa_i l}$, $\mathcal{L}_0; \mathcal{L}_i \in \mathbb{R}^{m\times l}$, $\widehat{\mathcal{L}}_i \in \mathbb{R}^{m\times\kappa_i l}$, and matrices $Q_i; R_i \in \mathbb{S}^n$, $i \in \mathbb{N}_\nu$ and matrices $P_1 \in \mathbb{S}^n$, $P_2 \in \mathbb{R}^{n\times dn}$, $P_3 \in \mathbb{S}^{dn}$ with $d = \sum_{i=1}^{\nu} d_i$ such that

$$\begin{bmatrix}P_1 & P_2\\ * & P_3 + \mathsf{diag}_{i=1}^{\nu}(I_{d_i} \otimes Q_i)\end{bmatrix} \succ 0, \qquad (31a)$$

$$\mathbf{Q} = \operatorname*{diag}_{i=1}^{\nu} Q_i \succ 0, \quad \mathbf{R} = \operatorname*{diag}_{i=1}^{\nu} R_i \succ 0, \qquad (31b)$$

$$\begin{bmatrix}\boldsymbol{\Psi} & \boldsymbol{\Sigma}^{\top}\widetilde{J}^{\top}\\ * & J_1\end{bmatrix} = \mathsf{Sy}\left[\mathbf{P}^{\top}\boldsymbol{\Pi}\right] + \boldsymbol{\Phi} \prec 0, \qquad (31c)$$

where $\boldsymbol{\Sigma} = \mathbf{C} + \mathbf{L}_2\left[(I_\beta \otimes \mathcal{C}) \oplus \mathsf{O}_q\right]$ with $\mathbf{C}, \mathbf{L}_2$ in (26d)–(26f) and $\beta$ in (25), and

$$\boldsymbol{\Psi} = \mathsf{Sy}\left(S^{\top}\begin{bmatrix}P_1 & P_2\\ * & P_3\end{bmatrix}\begin{bmatrix} & \boldsymbol{\Omega}\\ \mathbf{M}\otimes I_n & \mathsf{O}_{dn,(\mu n+q)}\end{bmatrix}\right. \qquad (32a)$$
$$\left. - \begin{bmatrix}\mathsf{O}_{(\beta n),m}\\ J_2^{\top}\end{bmatrix}\boldsymbol{\Sigma}\right) + \Xi, \qquad (32b)$$

$$S = \begin{bmatrix} I_n & \mathsf{O}_{n,\nu n} & \mathsf{O}_{n,\varkappa n} & \mathsf{O}_{n,\mu n} & \mathsf{O}_{n,q}\\ \mathsf{O}_{dn,n} & \mathsf{O}_{dn,\nu n} & \widehat{I} & \mathsf{O}_{dn,\mu n} & \mathsf{O}_{dn,q}\end{bmatrix}, \qquad (32c)$$

$$\Xi = \Big[(\mathbf{Q} + \mathbf{R}\Lambda) \oplus \mathsf{O}_n \oplus \mathsf{O}_{\kappa n} \oplus \mathsf{O}_q\Big] - \cdots \qquad (32d)$$
$$\Big[\mathsf{O}_n \oplus \mathbf{Q} \oplus \Big[\operatorname*{diag}_{i=1}^{\nu}\left(I_{\varkappa_i} \otimes R_i\right)\Big] \oplus \Big[\operatorname*{diag}_{i=1}^{\nu}\left(I_{\mu_i} \otimes R_i\right)\Big] \oplus J_3\Big],$$

$$\widehat{I} = \operatorname*{diag}_{i=1}^{\nu}\sqrt{\mathfrak{F}_i^{-1}}\widetilde{I}_i\sqrt{\mathfrak{H}_i} \otimes I_n, \ \ \widetilde{I}_i = \begin{bmatrix}\mathsf{O}_{d_i,\delta_i} & I_{d_i}\end{bmatrix}, \qquad (32e)$$

$$\Lambda = \operatorname*{diag}_{i=1}^{\nu} \acute{r}_i I_n, \ \ \acute{r}_i = r_i - r_{i-1}, \qquad (32f)$$

$$\mathbf{M} = \Big[\mathsf{diag}_{i=1}^{\nu}\sqrt{\mathfrak{F}_i^{-1}}\boldsymbol{f}_i(-r_{i-1}) \quad \mathbf{0}_d \quad \mathsf{O}_{d,\varkappa}\Big] - \cdots$$
$$\Big[\mathbf{0}_d \quad \mathsf{diag}_{i=1}^{\nu}\sqrt{\mathfrak{F}_i^{-1}}\boldsymbol{f}_i(-r_i) \quad \mathsf{diag}_{i=1}^{\nu}\sqrt{\mathfrak{F}_i^{-1}}M_i\sqrt{\mathfrak{H}_i}\Big] \quad (32g)$$

with $\varkappa, \mu, \kappa$ in (23)–(25), and $\kappa_i, \varkappa_i, \mu_i, M_i$ in Proposition 1, and $\boldsymbol{\Omega} = \mathbf{A} + \mathbf{L}_1\left[(I_\beta \otimes \mathcal{C}) \oplus \mathsf{O}_q\right]$ with $\mathbf{A}, \mathbf{L}_1$ in (26c)–(26e), and $\mathfrak{F}_i := \int_{\mathcal{I}_i} \boldsymbol{f}_i(\tau)\boldsymbol{f}_i^{\top}(\tau)\mathsf{d}\tau$, $\forall i \in \mathbb{N}_\nu$. Moreover,

$$\mathbf{P} = \Big[P_1 \quad \mathsf{O}_{n,\nu n} \quad P_2\widehat{I} \quad \mathsf{O}_{n,(\mu n+q+m)}\Big], \ \boldsymbol{\Pi} = \begin{bmatrix}\boldsymbol{\Omega} & \mathsf{O}_{n,m}\end{bmatrix}, \ (33)$$

$$\mathbf{\Phi} = \mathsf{Sy}\left(\begin{bmatrix} P_2 \\ \mathsf{O}_{\nu n,dn} \\ \widehat{I}^\top P_3 \\ \mathsf{O}_{(\mu n+q+m),dn} \end{bmatrix}\begin{bmatrix}\mathbf{M}\otimes I_n & \mathsf{O}_{dn,(\mu n+q+m)}\end{bmatrix} + \begin{bmatrix} \mathsf{O}_{(\beta n),m} \\ -J_2^\top \\ \widetilde{J} \end{bmatrix}\begin{bmatrix}\mathbf{\Sigma} & \mathsf{O}_m\end{bmatrix}\right) + \Xi \oplus (-J_1)\,.$$

Finally, the number of unknowns in the constraints of this theorem is $(0.5d^2+0.5d+\nu+1)n^2+(0.5d+1+\beta l)n+\beta lm \in \mathcal{O}(d^2n^2)$. Note that $\mathsf{diag}_{i=1}^{\nu} X_i \otimes I_n$ stands for $[\mathsf{diag}_{i=1}^{\nu} X_i]\otimes I_n$ given the order of operation defined in the introduction.

*Remark 6:* Numerical solutions to (31a)–(31c) can always ensure stability and dissipativity for (17), as (25) is equivalent to (17). This is different from numerically solving operator Riccati equations [55], whose numerical solutions do not strictly guarantee the stability of the original CLS due to the use of finite-dimensional approximations.

The inequality in (31c) is bilinear due to the products between the estimator gains and $P_1, P_2$, which cannot be computed by standard SDP numerical solvers. To address this issue, we propose the following theorem, where (31c) is made convex using [40, Lem. 3.1] that could decouple the products between two decision variables [56], [57].

*Lemma 3 (Projection Lemma):* [40, Lem. 3.1] Given $n; p; q \in \mathbb{N}$, $\Pi \in \mathbb{S}^n$, $P \in \mathbb{R}^{q\times n}$, $Q \in \mathbb{R}^{p\times n}$, there exists $\Theta \in \mathbb{R}^{p\times q}$ such that the following two assertions are equivalent :

$$\Pi + P^\top \Theta^\top Q + Q^\top \Theta P \prec 0, \tag{34a}$$

$$P_\perp^\top \Pi P_\perp \prec 0 \;\text{ and }\; Q_\perp^\top \Pi Q_\perp \prec 0, \tag{34b}$$

where the columns of $P_\perp, Q_\perp$ contain some bases of the null space of $P, Q$, respectively, satisfying $PP_\perp = \mathsf{O}, QQ_\perp = \mathsf{O}$.

*Theorem 2:* Given Proposition 1 and (8) with known scalars $\{\alpha_i\}_{i=1}^{\beta} \subset \mathbb{R}$, then system (25) with SRF (30) is dissipative and the trivial solution to (25) with $\boldsymbol{w}(t) \equiv \mathbf{0}_q$ is exponentially stable for any $\boldsymbol{\psi}(\cdot) \in \mathcal{C}(\mathcal{J};\mathbb{R}^n)$ in (25) if there exists $P_1; W \in \mathbb{S}^n$, $P_2 \in \mathbb{R}^{n\times dn}$, $P_3 \in \mathbb{S}^{dn}$ and $Q_i; R_i \in \mathbb{S}^n$, $U_0; U_i \in \mathbb{R}^{n\times l}$, $\widehat{U}_i \in \mathbb{R}^{n\times \kappa_i l}$, $\mathcal{L}_0; \mathcal{L}_i \in \mathbb{R}^{m\times l}$, $\widehat{\mathcal{L}}_i \in \mathbb{R}^{m\times \kappa_i l}$, $i \in \mathbb{N}_\nu$ such that the linear matrix inequalities in (31a)–(31b) and

$$\mathsf{Sy}\left(\begin{bmatrix} I_n \\ [\alpha_i I_n]_{i=1}^{\beta} \\ \mathsf{O}_{(q+m),n} \end{bmatrix}\begin{bmatrix}-W & \widehat{\mathbf{\Pi}}\end{bmatrix}\right) + \begin{bmatrix} \mathsf{O}_n & \mathbf{P} \\ * & \mathbf{\Phi} \end{bmatrix} \prec 0 \tag{35}$$

hold with $\mathbf{P}, \mathbf{\Phi}$ given by Theorem 1, where

$$\widehat{\mathbf{\Pi}} = \begin{bmatrix} W\mathbf{A} + \mathbf{U}\left[(I_\beta \otimes \mathcal{C}) \oplus \mathsf{O}_q\right] & \mathsf{O}_{n,m}\end{bmatrix}, \tag{36a}$$

$$\mathbf{U} = \Big[U_0 \quad [\![U_i]\!]_{i=1}^{\nu} \quad \left[\!\!\left[\widehat{U}_i\left(T_i \otimes I_l\right)\right]\!\!\right]_{i=1}^{\nu} \cdots \cdots \left[\!\!\left[\widehat{U}_i\left(\widetilde{T}_i \otimes I_l\right)\right]\!\!\right]_{i=1}^{\nu} \quad \mathsf{O}_{n,q}\Big] \tag{36b}$$

with $\mathbf{A}$ in (26c) and $T_i, \widetilde{T}_i$ in (20)–(21). The state estimator gains are obtained as $L_i = W^{-1}U_i$, $\widehat{L}_i = W^{-1}\widehat{U}_i$. Finally, the number of unknowns is $(0.5d^2+0.5d+\nu+1)n^2+(0.5d+1+\nu+\beta l)n+\beta ml \in \mathcal{O}(d^2n^2)$ with $\beta$ in (25).

*Remark 7:* For scalars $\{\alpha_i\}_{i=1}^{\beta} \subset \mathbb{R}$, we can let $\alpha_i = 0$ for $i = 2, \ldots, \beta$ with a tunable parameter $\alpha_1 \in \mathbb{R}\setminus\{0\}$. Note that $\alpha_1 \neq 0$ is necessary as $\alpha_1 = 0$ may render the $A_0$ related-diagonal-block in (35) infeasible.

## A. Sequential Convex SDP for solving BMI

Although Theorem 2 furnishes convex solutions to the DSE design problem, step (66) can introduce conservatism. Hence, it is preferable to find methods that can compute (31c) directly. We propose a sequential convex SDP algorithm as an application of the concept of inner convex approximations [58], which ensures monotonic convergence to a local optimum given the general conclusions in [58, Th. 3.2]. Each iteration in our algorithm solves a convex SDP problem, and the initial point of the algorithm can be provided by a feasible solution to (31a)–(31b) and (35), maximizing the utility of Theorem 2.

First, we can reformulate (31c) as

$$\mathsf{Sy}\left[\mathbf{P}^\top \mathbf{\Pi}\right] + \mathbf{\Phi} = \mathsf{Sy}\left(\mathbf{P}^\top \mathbf{L}\left[(I_\beta \otimes \mathcal{C}) \oplus \mathsf{O}_{q+m}\right]\right) + \widehat{\mathbf{\Phi}} \prec 0, \tag{37}$$

where $\mathbf{L} := \begin{bmatrix}\mathbf{L}_1 & \mathsf{O}_{n,m}\end{bmatrix}$ and $\widehat{\mathbf{\Phi}} := \mathsf{Sy}\left(\mathbf{P}^\top \begin{bmatrix}\mathbf{A} & \mathsf{O}_{n,m}\end{bmatrix}\right) + \mathbf{\Phi}$ with $\mathbf{P}$ in (33) and $\mathbf{A}, \mathbf{L}_1$ in (26c)–(26e). Note that $\widehat{\mathbf{\Phi}}$ contains no nonconvex terms. Utilizing the conclusion of [58, Example 3], the function $\Delta\big(\bullet, \widetilde{\mathbf{G}}, \bullet, \widetilde{\mathbf{N}}\big)$ denoted by

$$\mathbb{S}^{\ell} \ni \Delta\big(\mathbf{G}, \widetilde{\mathbf{G}}, \mathbf{N}, \widetilde{\mathbf{N}}\big) := [*]\left[Z \oplus (I_n - Z)\right]^{-1}\begin{bmatrix}\mathbf{G}-\widetilde{\mathbf{G}} \\ \mathbf{N}-\widetilde{\mathbf{N}}\end{bmatrix} + \mathsf{Sy}\left(\widetilde{\mathbf{G}}^\top \mathbf{N} + \mathbf{G}^\top \widetilde{\mathbf{N}} - \widetilde{\mathbf{G}}^\top \widetilde{\mathbf{N}}\right) + \mathbf{T} \tag{38a}$$

with $Z \oplus (I_n - Z) \succ 0$ satisfying

$$\forall \mathbf{G}; \widetilde{\mathbf{G}} \in \mathbb{R}^{n\times \ell},\ \forall \mathbf{N}; \widetilde{\mathbf{N}} \in \mathbb{R}^{n\times \ell},\quad \mathbf{T} + \mathsf{Sy}\left(\mathbf{G}^\top \mathbf{N}\right) = \Delta\left(\mathbf{G}, \mathbf{G}, \mathbf{N}, \mathbf{N}\right) \preceq \Delta\big(\mathbf{G}, \widetilde{\mathbf{G}}, \mathbf{N}, \widetilde{\mathbf{N}}\big), \tag{38b}$$

is a PSD-convex mapping [59, Eq. (6)] and PSD-convex overestimate [58, Definition 2.2] of $\Delta(\mathbf{G}, \mathbf{G}, \mathbf{N}, \mathbf{N}) = \mathbf{T} + \mathsf{Sy}\left[\mathbf{G}^\top \mathbf{N}\right]$ w.r.t the parameterization

$$\mathsf{Col}\left[\mathbf{vec}(\widetilde{\mathbf{G}}), \mathbf{vec}(\widetilde{\mathbf{N}})\right] = \mathsf{Col}\left[\mathbf{vec}(\mathbf{G}), \mathbf{vec}(\mathbf{N})\right].$$

Now let $\mathbf{T} = \widehat{\mathbf{\Phi}}$, $\mathbf{G} = \mathbf{P}$ and $\widetilde{P}_1 \in \mathbb{S}^n$, $\widetilde{P}_2 \in \mathbb{R}^{n\times dn}$ and

$$\begin{aligned}
&\widetilde{\mathbf{G}} = \widetilde{\mathbf{P}} = \begin{bmatrix}\widetilde{P}_1 & \mathsf{O}_{n,\nu n} & \widetilde{P}_2\widehat{I} & \mathsf{O}_{n,(\mu n+q+m)}\end{bmatrix},\\
&\mathbf{N} = \mathbf{L}\left[(I_\beta \otimes \mathcal{C}) \oplus \mathsf{O}_{q+m}\right],\quad \widetilde{\mathbf{N}} = \widetilde{\mathbf{L}}\left[(I_\beta \otimes \mathcal{C}) \oplus \mathsf{O}_{q+m}\right],\\
&\widetilde{\mathbf{L}} = \Big[V_0 \quad [\![V_i]\!]_{i=1}^{\nu} \quad \left[\!\!\left[\widehat{V}_i\left(T_i \otimes I_l\right)\right]\!\!\right]_{i=1}^{\nu} \cdots \cdots \left[\!\!\left[\widehat{V}_i\left(\widetilde{T}_i \otimes I_l\right)\right]\!\!\right]_{i=1}^{\nu} \quad \mathsf{O}_{n,(q+m)}\Big],\\
&\mathfrak{L} = \Big[[\![L_i]\!]_{i=0}^{\nu} \quad \left[\!\!\left[\widehat{L}_i\right]\!\!\right]_{i=1}^{\nu}\Big],\ \widetilde{\mathfrak{L}} = \Big[[\![V_i]\!]_{i=0}^{\nu} \quad \left[\!\!\left[\widehat{V}_i\right]\!\!\right]_{i=1}^{\nu}\Big]
\end{aligned} \tag{39}$$

for (38a) with $\ell = \beta n + q + m$ and $T_i, \widetilde{T}_i$ in (20)–(21) and $\widehat{\mathbf{\Phi}}, \mathbf{L}$ in (37), where $V_i \in \mathbb{R}^{n\times l}$, $\widehat{V}_i \in \mathbb{R}^{n\times \kappa_i l}$. Then we get

$$\begin{aligned}
&\mathcal{S}\left(\mathbf{P}, \mathbf{P}, \mathfrak{L}, \mathfrak{L}\right) = \widehat{\mathbf{\Phi}} + \mathsf{Sy}\left[\mathbf{P}^\top \mathbf{L}\left[(I_\beta \otimes \mathcal{C}) \oplus \mathsf{O}_{p+m}\right]\right] \text{ in (37)}\\
&\preceq \mathcal{S}\left(\mathbf{P}, \widetilde{\mathbf{P}}, \mathfrak{L}, \widetilde{\mathfrak{L}}\right) := \widehat{\mathbf{\Phi}} + \mathsf{Sy}\left(\widetilde{\mathbf{P}}^\top \mathbf{N} + \mathbf{P}^\top \widetilde{\mathbf{N}} - \widetilde{\mathbf{P}}^\top \widetilde{\mathbf{N}}\right)\\
&\quad + \begin{bmatrix}\mathbf{P}^\top - \widetilde{\mathbf{P}}^\top & \mathbf{N}^\top - \widetilde{\mathbf{N}}^\top\end{bmatrix}\left[Z \oplus (I_n - Z)\right]^{-1}[*]
\end{aligned} \tag{40}$$

by (38b), where $\mathcal{S}(\bullet, \widetilde{\mathbf{P}}, \bullet, \widetilde{\mathfrak{L}})$ is a PSD-convex overestimate of $\mathcal{S}\left(\mathbf{P}, \mathbf{P}, \mathfrak{L}, \mathfrak{L}\right)$ in (37) w.r.t the parameterization

$$\mathsf{Col}\left[\mathbf{vec}(\widetilde{\mathbf{P}}), \mathbf{vec}(\widetilde{\mathfrak{L}})\right] = \mathsf{Col}\left[\mathbf{vec}(\mathbf{P}), \mathbf{vec}(\mathfrak{L})\right].$$

With (40), we can infer (37) from $\mathcal{S}\left(\mathbf{P}, \widetilde{\mathbf{P}}, \mathfrak{L}, \widetilde{\mathfrak{L}}\right) \prec 0$ by $Z \oplus (I_n - Z) \succ 0$. Finally, we have $\mathcal{S}\left(\mathbf{P}, \widetilde{\mathbf{P}}, \mathfrak{L}, \widetilde{\mathfrak{L}}\right) \prec 0 \iff$

$$\begin{bmatrix} \widehat{\mathbf{\Phi}} + \mathsf{Sy}\left(\widetilde{\mathbf{P}}^\top \mathbf{N} + \mathbf{P}^\top \widetilde{\mathbf{N}} - \widetilde{\mathbf{P}}^\top \widetilde{\mathbf{N}}\right) & \mathbf{P}^\top - \widetilde{\mathbf{P}}^\top & \mathbf{N}^\top - \widetilde{\mathbf{N}}^\top \\ * & -Z & \mathsf{O}_n \\ * & * & Z - I_n \end{bmatrix} \prec 0 \quad (41)$$

with $\mathbf{N}$, $\widetilde{\mathbf{N}}$ in (39) by a use of the Schur complement with $Z \oplus (I_n - Z) \succ 0$, where (41) implies (37). Note that the constraint in (41) is convex if $\widetilde{\mathbf{P}}, \widetilde{\mathfrak{L}}$ are known.

By synthesizing the above procedures in accordance with the expositions in [58], we formulate Algorithm 1, where $\mathbf{x}$ comprises all the variables in Theorem 1 and $Z$ in (41), and $\|\mathbf{x}\|_\infty := \max_i x_i$. Finally, $\rho_1, \rho_2$ and $\varepsilon$ are given values for regularizations and for indicating error tolerance, respectively.

**Algorithm 1:** Sequential Convex SDP for Theorem 1

**begin**
- **utilize** Theorem 2 **compute** $\mathfrak{L}$
- **utilize** Theorem 1 with $\mathfrak{L}$ **compute** $\mathbf{Y} = \begin{bmatrix} P_1 & P_2 \end{bmatrix}$
- **utilize** Theorem 1 with $\mathbf{Y}$ **compute** $\mathfrak{L}$ again.
- **update** $\widetilde{\mathbf{Y}} = \begin{bmatrix} \widetilde{P}_1 & \widetilde{P}_2 \end{bmatrix} \longleftarrow \mathbf{Y} \quad \widetilde{\mathfrak{L}} \longleftarrow \mathfrak{L}$,
- **solve** $\mathcal{P} : \min_{\mathbf{x}} \mathsf{tr}\left[\rho_1 [*]\left(\mathbf{Y} - \widetilde{\mathbf{Y}}\right)\right] + \mathsf{tr}\left[\rho_2 [*]\left(\mathfrak{L} - \widetilde{\mathfrak{L}}\right)\right]$ subject to (31a)–(31b), (41) with (39) and the parameters in Theorem 1, **return** $\mathbf{Y}$ and $\mathfrak{L}$
- **while** $\dfrac{\left\| \begin{bmatrix} \mathbf{vec}(\mathbf{Y}) \\ \mathbf{vec}(\mathfrak{L}) \end{bmatrix} - \begin{bmatrix} \mathbf{vec}(\widetilde{\mathbf{Y}}) \\ \mathbf{vec}(\widetilde{\mathfrak{L}}) \end{bmatrix} \right\|_\infty}{\left\| \begin{bmatrix} \mathbf{vec}(\widetilde{\mathbf{Y}}) \\ \mathbf{vec}(\widetilde{\mathfrak{L}}) \end{bmatrix} \right\|_\infty + 1} \geq \varepsilon$ **do**
  - **update** $\widetilde{\mathbf{Y}} \longleftarrow \mathbf{Y}, \quad \widetilde{\mathfrak{L}} \longleftarrow \mathfrak{L}$,
  - **solve** $\mathcal{P}$ again and **return** $\mathbf{Y}$ and $\mathfrak{L}$
- **end**

**end**

The following diagram demonstrates the relations between the proposed theorems and iterative algorithm. Together, they can be used as a package to solve the DSE problem for (25).

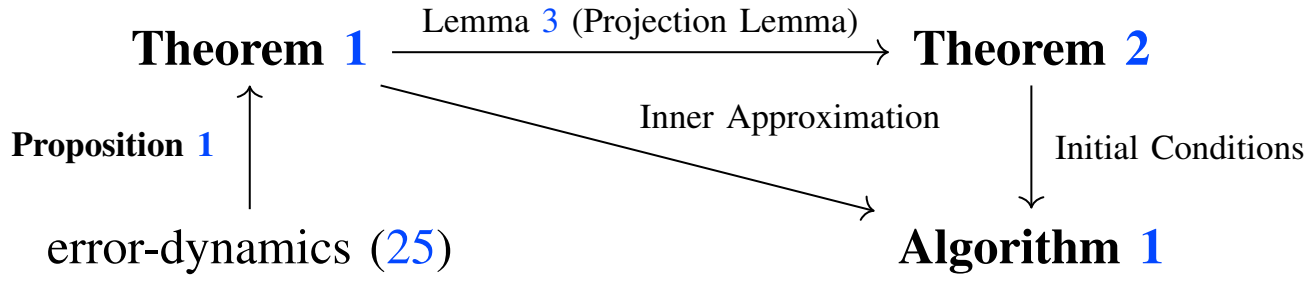

Schematic Diagram of DSE design for (25)

## IV. Numerical Examples and Experiments

Two numerical examples were tested to show the effectiveness of our framework. All computations were carried out in MATLAB© using Yalmip [60] as the optimization parser, and Mosek [24] as the numerical solver for SDP.

### A. DSE Design of an LTDS with Multiple Delays

Consider a system in the form of (1) with $r_1 = 1$, $r_2 = 1.7$ and $\mathcal{I}_1 = [-1, 0]$, $\mathcal{I}_2 = [-1.7, -1]$ and state space matrices

$$A_0 = \begin{bmatrix} -3 & 0.4 \\ 0 & 0.2 \end{bmatrix}, \ A_1 = \begin{bmatrix} 0.2 & 0.8 \\ -0.4 & -0.5 \end{bmatrix}, \ A_2 = \begin{bmatrix} -0.2 & 0.1 \\ 0.3 & 0.2 \end{bmatrix}$$

$$\widetilde{A}_1(\tau) = \begin{bmatrix} 0.1 + 3\sin 17\tau & 0.8\mathsf{e}^{\sin 17\tau} - 0.3\mathsf{e}^{\cos 17\tau} \\ 0.3 + \sin\frac{1}{\tau - 0.1} + 0.5 & 3\sin(17\tau) + \varsigma(t) \end{bmatrix}$$

$$\widetilde{A}_2(\tau) = \begin{bmatrix} -10\cos 21\tau & 0.3\mathsf{e}^{\cos 21\tau} - \cos\frac{1}{\tau + 0.9} - 0.5 \\ 0.1\mathsf{e}^{\sin 21\tau} & 0.2 - 10\cos 21\tau \end{bmatrix}$$

$$C_0 = \begin{bmatrix} 0 & 1 \end{bmatrix}, \ C_1 = \begin{bmatrix} -0.1 & 0 \end{bmatrix}, \ C_2 = \begin{bmatrix} 0 & 0.1 \end{bmatrix} \quad (42)$$

$$\widetilde{C}_1(\tau) = \begin{bmatrix} 0.1\mathsf{e}^{\sin 17\tau} + 0.1\mathsf{e}^{\cos 17\tau} & 0.4\left(\sin\frac{1}{\tau - 0.1} + 0.5\right) + 1 \end{bmatrix}$$

$$\widetilde{C}_2(\tau) = \begin{bmatrix} 0.2 & \mathsf{e}^{\cos 21\tau} + 0.2\mathsf{e}^{\sin 21\tau} + 0.3 \end{bmatrix}, \ \mathcal{C} = \begin{bmatrix} 0 & 1 \end{bmatrix}$$

$$D_1 = \begin{bmatrix} 1.5 \\ 0.5 \end{bmatrix}, \ D_2 = 0.7, \ D_3 = \begin{bmatrix} 0.2 \\ 0.3 \end{bmatrix}, \ D_4 = 0.2$$

with $\mathfrak{f}_i = \mathbf{0}_2$ and $n = 2, m = l = q = 1$, where $\widetilde{\forall} t \in \mathbb{R}$, $\varsigma(t) = n(t)\mathbb{1}_{\mathbb{Q}}(t) = 0$. Signal $\varsigma(t)$ may stand for glitches or other anomalies in the DD kernel matrix $\widetilde{A}_1(\tau)$, where $n(t)$ can be a function or stochastic process and $\mathbb{1}_{\mathbb{Q}}(t)$ is the indicator function over rational number set $\mathbb{Q}$. The incorporation of $\varsigma(t)$ is not achievable if we were to use the traditional derivative for $\dot{\boldsymbol{x}}(t)$ in (1), even if we were to consider the time-varying FDE framework in [53], [61]. This is because the right hand side of $\dot{\boldsymbol{x}}(t)$ in (1) in this case is neither continuous nor piecewise continuous in $t$, but satisfying the Carathéodory condition [12, page 58] as an FDE in the extended sense. Since set $\mathbb{Q}$ is countable, signal $\varsigma(t) = n(t)\mathbb{1}_{\mathbb{Q}}(t)$ is regarded as zero w.r.t the Lebesgue measure. In fact, $\varsigma(t)$ can be added to any matrix in (42) and be treated as 0 by the proposed approach. This serves as an excellent example that can show the benefits of employing the Carathéodory formulation in the modeling of TDSs. By utilizing the numerical toolbox https://cdlab.uniud.it/software#eigAM-eigTMN of the spectral method in [62], we found that (1) with (42) and $\boldsymbol{w}(t) \equiv \mathbf{0}_q$ is unstable. The parameters in (42) were intentionally chosen to show that the validity of our estimator is not affected by the stability of (1). Moreover, we utilize the $\mathcal{L}^2$ gain performance

$$\gamma > 0, \quad J_1 = -\gamma I_2, \quad \widetilde{J} = I_2, \quad J_2 = \mathbf{0}_2, \quad J_3 = \gamma \quad (43)$$

for (30) with $\gamma > 0$ as the variable to be minimized.

*Remark 8:* The parameters in (42) were selected with a sufficient degree of complexity to illustrate the strength of the proposed method. Note that our approach can handle practical examples listed in Remark 1, whose DD matrices are **simpler** than those in (42). To the best of our knowledge, no existing methods may effectively compute solutions for the DSE design problem of (1) with (42) due to the complexity of the DD matrices with non-commensurate delays and non-Hurwitz $A_0$.

Our goal is to compute the gains of the estimator in (12) to stabilize the error-dynamics in (25) while minimizing $\gamma$. Observing the DD kernels in (42), let $\sigma_1; \sigma_2 \in \mathbb{N}_0$, $\lambda_1; \lambda_2 \in \mathbb{N}$, and $\boldsymbol{\varphi}_1(\tau) = \sin\frac{1}{\tau - 0.1} + 0.5$, $\boldsymbol{\varphi}_2(\tau) = \cos\frac{1}{\tau + 0.9} + 0.5$ and

$$\boldsymbol{f}_1(\tau) = \begin{bmatrix} \left[\tau^i\right]_{i=0}^{\sigma_1} \\ \left[\sin 17i\tau\right]_{i=1}^{\lambda_1} \\ \left[\cos 17i\tau\right]_{i=1}^{\lambda_1} \end{bmatrix}, \quad \boldsymbol{f}_2(\tau) = \begin{bmatrix} \left[\tau^i\right]_{i=0}^{\sigma_2} \\ \left[\sin 21i\tau\right]_{i=1}^{\lambda_2} \\ \left[\cos 21i\tau\right]_{i=1}^{\lambda_2} \end{bmatrix} \quad (44a)$$

$$\boldsymbol{\phi}_1(\tau) = \begin{bmatrix} \mathsf{e}^{\sin 17\tau} \\ \mathsf{e}^{\cos 17\tau} \end{bmatrix}, \qquad \boldsymbol{\phi}_2(\tau) = \begin{bmatrix} \mathsf{e}^{\sin 21\tau} \\ \mathsf{e}^{\cos 21\tau} \end{bmatrix}$$

for (4)–(6) with $d_i = 2\lambda_i + \sigma_i + 1$, $\mu_i = 2$, $\delta_i = 1$ and

$$M_1 = \begin{bmatrix} \mathbf{0}_{d_1} & \begin{bmatrix} \mathbf{0}_{\sigma_1}^\top & 0 \\ \mathsf{diag}_{i=1}^{\sigma_1} i & \mathbf{0}_{\sigma_1} \end{bmatrix} \oplus \begin{bmatrix} \mathsf{O}_{\lambda_1} & \mathsf{diag}_{i=1}^{\lambda_1} 17i \\ -\mathsf{diag}_{i=1}^{\lambda_1} 17i & \mathsf{O}_{\lambda_1} \end{bmatrix} \end{bmatrix}$$

$$M_2 = \begin{bmatrix} \mathbf{0}_{d_2} & \begin{bmatrix} \mathbf{0}_{\sigma_2}^\top & 0 \\ \mathsf{diag}_{i=1}^{\sigma_2} i & \mathbf{0}_{\sigma_2} \end{bmatrix} \oplus \begin{bmatrix} \mathsf{O}_{\lambda_2} & \mathsf{diag}_{i=1}^{\lambda_2} 21i \\ -\mathsf{diag}_{i=1}^{\lambda_2} 21i & \mathsf{O}_{\lambda_2} \end{bmatrix} \end{bmatrix}$$

for $\frac{\mathsf{d}\boldsymbol{f}_i(\tau)}{\mathsf{d}\tau} = M_i \boldsymbol{h}_i(\tau)$ in (5). Consequently, we can construct

$$\widehat{A}_1 = \begin{bmatrix} 0 & 0.8 & 0 & -0.3 & 0 & 0 & 0.1 & 0 & \mathbf{0}_{2\sigma_1}^\top & 3 & 0 & \mathbf{0}_{4\lambda_1-2}^\top \\ 0 & 0 & 0 & 0 & 1 & 0 & 0.3 & 0 & \mathbf{0}_{2\sigma_1}^\top & 0 & 3 & \mathbf{0}_{4\lambda_1-2}^\top \end{bmatrix}$$

$$\widehat{A}_2 = \begin{bmatrix} 0 & 0 & 0 & 0.3 & 0 & -1 & 0 & 0 & \mathbf{0}_{2\sigma_2+2\lambda_2}^\top & -10 & 0 & \mathbf{0}_{2\lambda_2-2}^\top \\ 0.1 & 0 & 0 & 0 & 0 & 0 & 0 & 0.2 & \mathbf{0}_{2\sigma_2+2\lambda_2}^\top & 0 & -10 & \mathbf{0}_{2\lambda_2-2}^\top \end{bmatrix} \quad (44b)$$

$$\widehat{C}_1 = \begin{bmatrix} 0.1 & 0 & 0.1 & 0 & 0 & 0.4 & 0 & 1 & \mathbf{0}_{2\sigma_2+4\lambda_2}^\top \end{bmatrix}$$

$$\widehat{C}_2 = \begin{bmatrix} 0 & 0.2 & 0 & 1 & 0 & 0 & 0.2 & 0.3 & \mathbf{0}_{2\sigma_2+4\lambda_2}^\top \end{bmatrix}$$

to satisfy the KSD in (4)–(6).

*Remark 9:* We selected the functions in (44a) for several reasons. The functions in $\boldsymbol{\phi}_i(\cdot)$ can be easily approximated by the use of trigonometric functions combined with polynomials, some of which exist in the DDs in (42). On the other hand, the functions in $\boldsymbol{\varphi}_i(\cdot)$ are directly pulled out, as it is difficult to approximate $\sin[1/(\tau - 0.1)]$ and $\cos[1/(\tau + 0.9)]$ with $\sigma_i, \lambda_i$ of reasonable values. The choice for (44a) involves all the ingredient of the proposed KSD concept, which balances both feasibility and the induced computational complexity.

Utilizing the parameters specified in (42)–(44b), we first apply Theorem 2 to (25) with $\sigma_1 = \sigma_2 = \lambda_1 = \lambda_2 = 1$ and $\alpha_i = 0, i = 2, \ldots, \beta$, $\alpha_1 = 30$. The program yields feasible estimator gains ensuring $\min \gamma = 14.2818$, whose values were utilized as an initial value for Algorithm 1. Note that the parameters in (8)–(10) were computed off-line via `vpaintegral` in MATLAB© to ensure numerical accuracy.

After running the remaining steps in Algorithm 1 with different sets of $\sigma_i, \lambda_i$, the results are summarized in Tables I and II, where SA stands for the spectral abscissa [62] of the error-dynamics with $\boldsymbol{w}(t) \equiv \mathbf{0}_q$, and NoIs denotes the number of iterations in the while loop. To show the advantages of the delay structures in (12), we also computed (See Tables III and IV) $\min \gamma$ by using a delay-free estimator with parameters $L_0, \mathcal{L}_0$ and $L_i = \widetilde{L}_i(\tau) = \mathsf{O}_{n,l}$, $\mathcal{L}_i = \widetilde{\mathcal{L}}_i(\tau) = \mathsf{O}_{m,l}$, $i \in \mathbb{N}_\nu$. The numerical results in Tables I–IV indicate that adding new WDLIFs (larger $\lambda_1, \lambda_2$) to $\boldsymbol{f}_i(\cdot)$ can enhance the feasibility of the synthesis conditions, leading to smaller $\min \gamma$, while the delay structures in (12) are also beneficial to decrease $\min \gamma$. Furthermore, we see that Algorithm 1 can produce estimator gains with significantly better performance than using Theorem 2 alone ($\min \gamma = 14.2818$).

TABLE I: ESTIMATOR PERFORMANCE COMPUTED WITH $\sigma_1 = \sigma_2 = \lambda_1 = \lambda_2 = 1$

| $\min \gamma$ | 0.5737 | 0.5656 | 0.5618 | 0.56 |
|---|---|---|---|---|
| SA | −0.5451 | −0.674 | −0.4782 | −0.5351 |
| NoIs | 0 | 5 | 10 | 15 |

TABLE II: ESTIMATOR PERFORMANCE COMPUTED WITH $\sigma_1 = \sigma_2 = \lambda_1 = \lambda_2 = 2$

| $\min \gamma$ | 0.5093 | 0.5048 | 0.5033 | 0.5024 |
|---|---|---|---|---|
| SA | −0.3706 | −0.4465 | −0.5561 | −0.5457 |
| NoIs | 0 | 5 | 10 | 15 |

TABLE III: ESTIMATOR PERFORMANCE COMPUTED WITH $L_i = \widetilde{L}_i(\tau) = \mathsf{O}_{n,l}$, $L_i = \widetilde{L}_i(\tau) = \mathsf{O}_{m,l}$ $\sigma_1 = \sigma_2 = \lambda_1 = \lambda_2 = 1$.

| $\min \gamma$ | 0.6461 | 0.6014 | 0.5830 | 0.5744 |
|---|---|---|---|---|
| SA | −0.9487 | −1.0469 | −1.1082 | −1.1138 |
| NoIs | 0 | 5 | 10 | 15 |

TABLE IV: ESTIMATOR PERFORMANCE COMPUTED WITH $L_i = \widetilde{L}_i(\tau) = \mathsf{O}_{n,l}$, $L_i = \widetilde{L}_i(\tau) = \mathsf{O}_{m,l}$ $\sigma_1 = \sigma_2 = \lambda_1 = \lambda_2 = 2$.

| $\min \gamma$ | 0.7641 | 0.5880 | 0.5296 | 0.5084 |
|---|---|---|---|---|
| SA | −0.5832 | −0.5305 | −0.4534 | −0.6122 |
| NoIs | 0 | 5 | 10 | 15 |

For simulations, we utilize the estimator gains with $\min \gamma = 0.5024$, NoIs = 15 in Table II. Assume $t_0 = 0$, $\boldsymbol{\zeta}(t) = \mathbf{0}_2, t < 0$, and $\widehat{\boldsymbol{\psi}}(\tau) = \begin{bmatrix} 2 & 1.8 \end{bmatrix}^\top, \boldsymbol{\psi}(\tau) = \begin{bmatrix} 0.5 & 1 \end{bmatrix}^\top$ as the initial conditions for (1) and (25), respectively, and $w(t) = \sin 20\pi(\mathit{1}(t) - \mathit{1}(t-3))$ as the disturbance where $\mathit{1}(t)$ is the Heaviside step function. To compute integrals numerically, we discretized all the DDs using the trapezoidal rule with 200 sampling points. We employ the Band-Limited White Noise block in Simulink to generate a white noise signal for $\varsigma(t)$ in (42) with `Sample time` $= 0.002s$ and the default values for `Seed`. As $\varsigma(t)$ can only be implemented as a discrete sequence $\varsigma(t) = \varsigma(kT)$ in the numerical simulation environment, hence $\varsigma(t)$ has only a **countable** number of nonzero data, which satisfies $\widetilde{\forall} t \geq 0$, $\varsigma(t) = \varsigma(kT) = 0$ as in (42). Our simulation was conducted in Simulink using `ode8` in fixed step size $0.002s$. To conclude, the plots in Figures 1–3 clearly show the validity of the resulting estimator, as both $\boldsymbol{e}(t)$ and $\boldsymbol{\zeta}(t)$ converge to zeros, even though the open-loop system is unstable and affected by a stochastic signal $\varsigma(t)$.

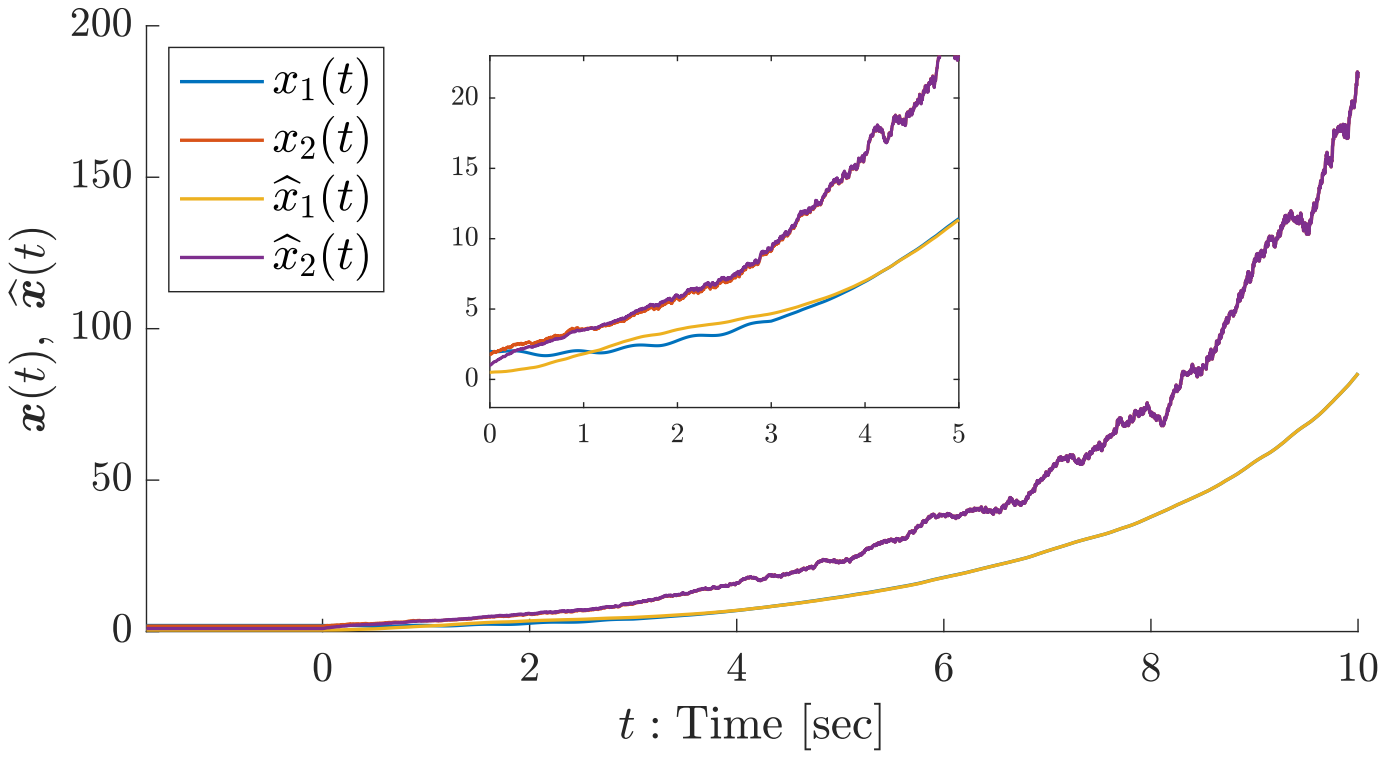

Fig. 1: Plots of $\boldsymbol{x}(t)$ and estimator state $\widehat{\boldsymbol{x}}(t)$

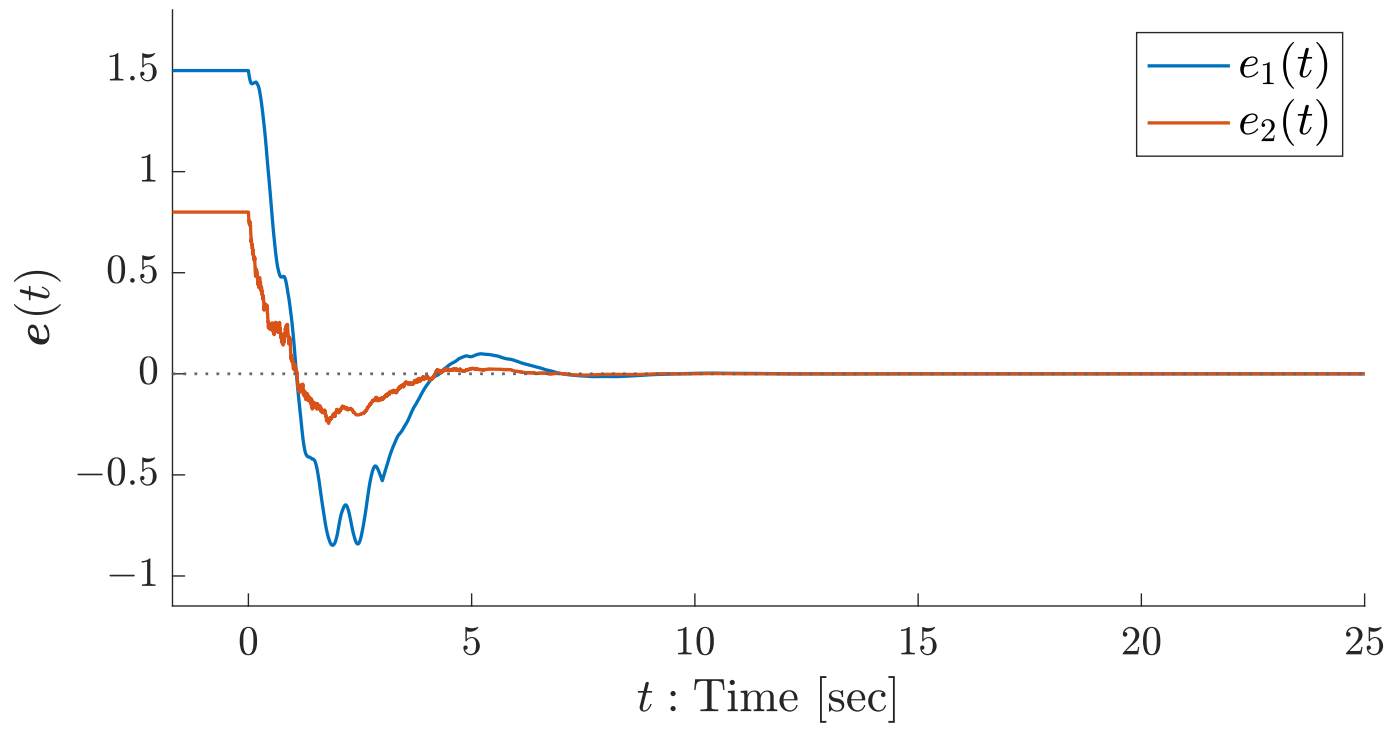

Fig. 2**:** Plots of error-dynamics $\boldsymbol{e}(t)$

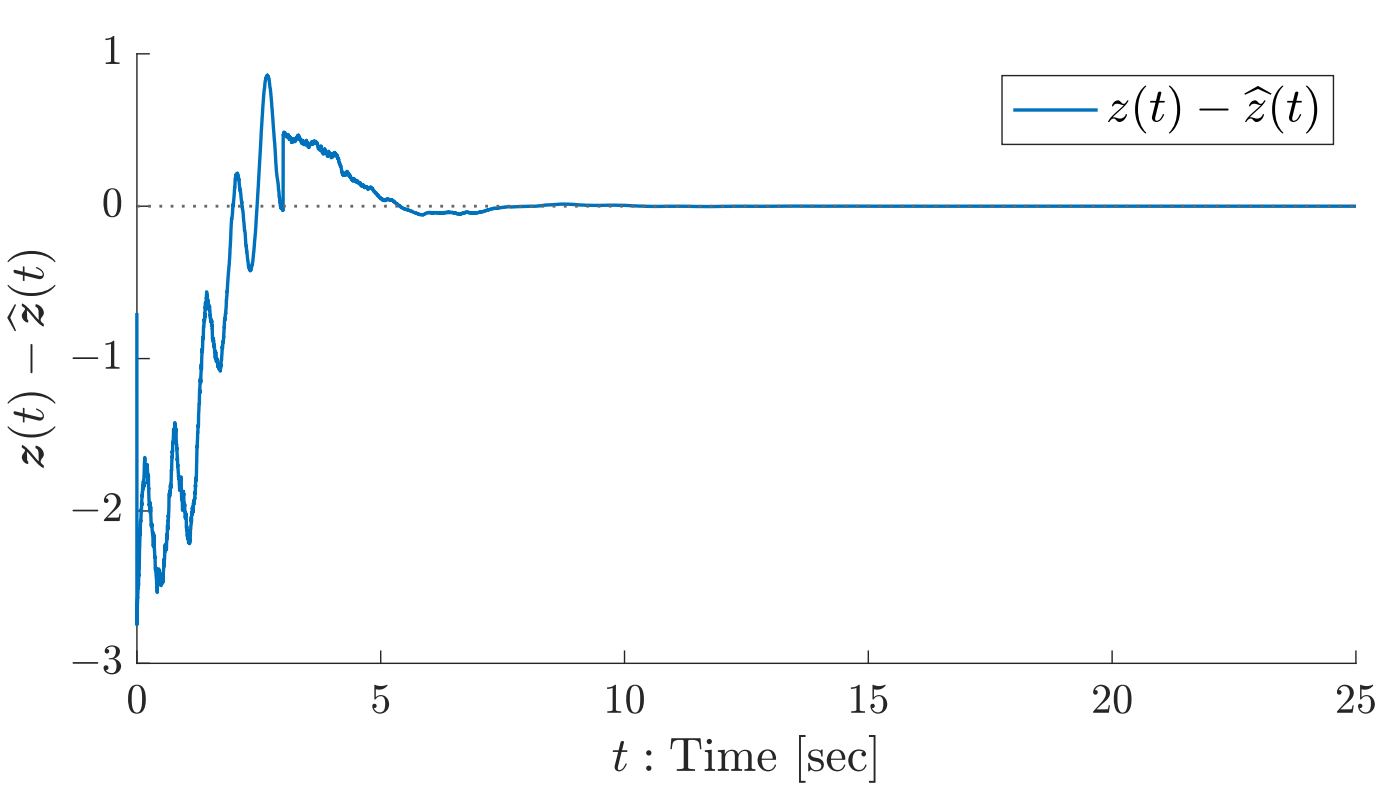

Fig. 4**:** Plots of output esitmation error $\boldsymbol{\zeta}(t) = \boldsymbol{z}(t) - \widehat{\boldsymbol{z}}(t)$

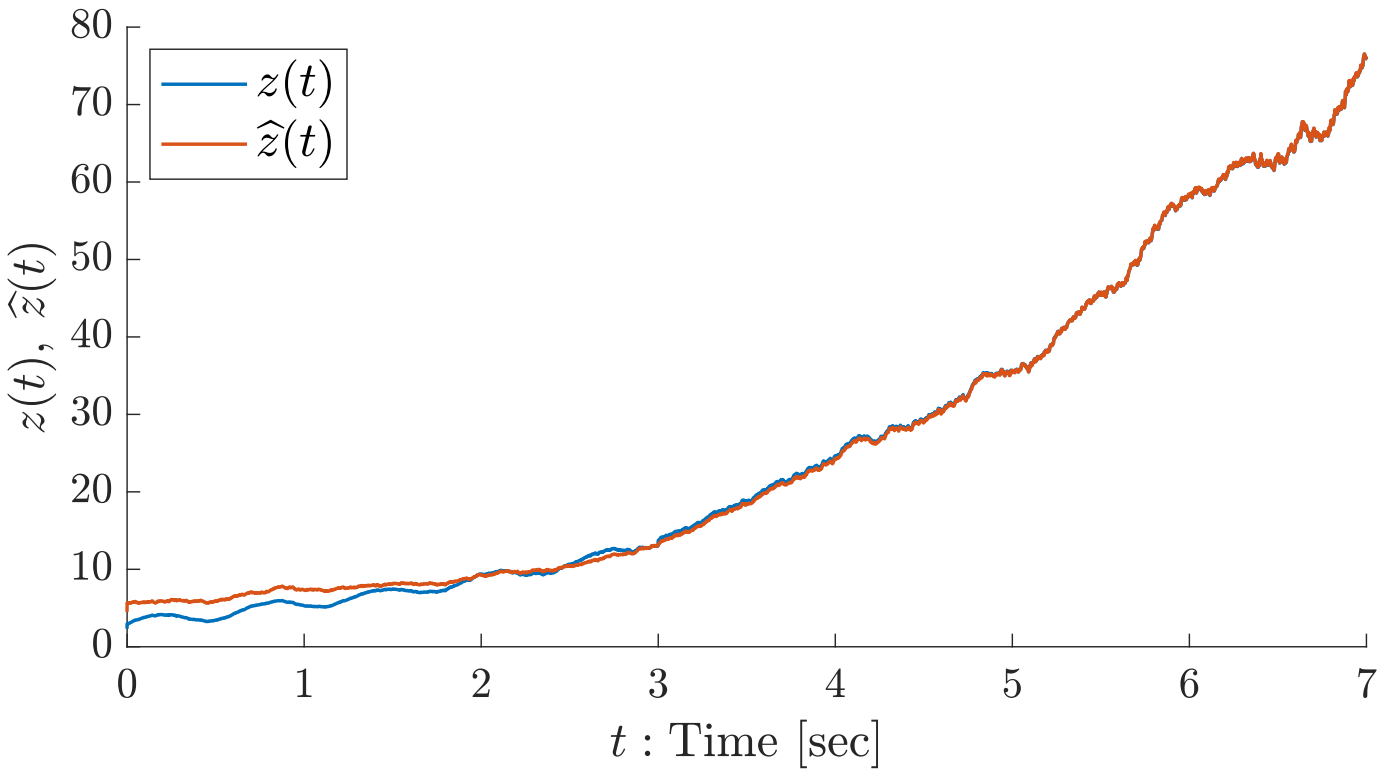

Fig. 3**:** Plots of regulated outputs $\boldsymbol{z}(t), \widehat{\boldsymbol{z}}(t)$

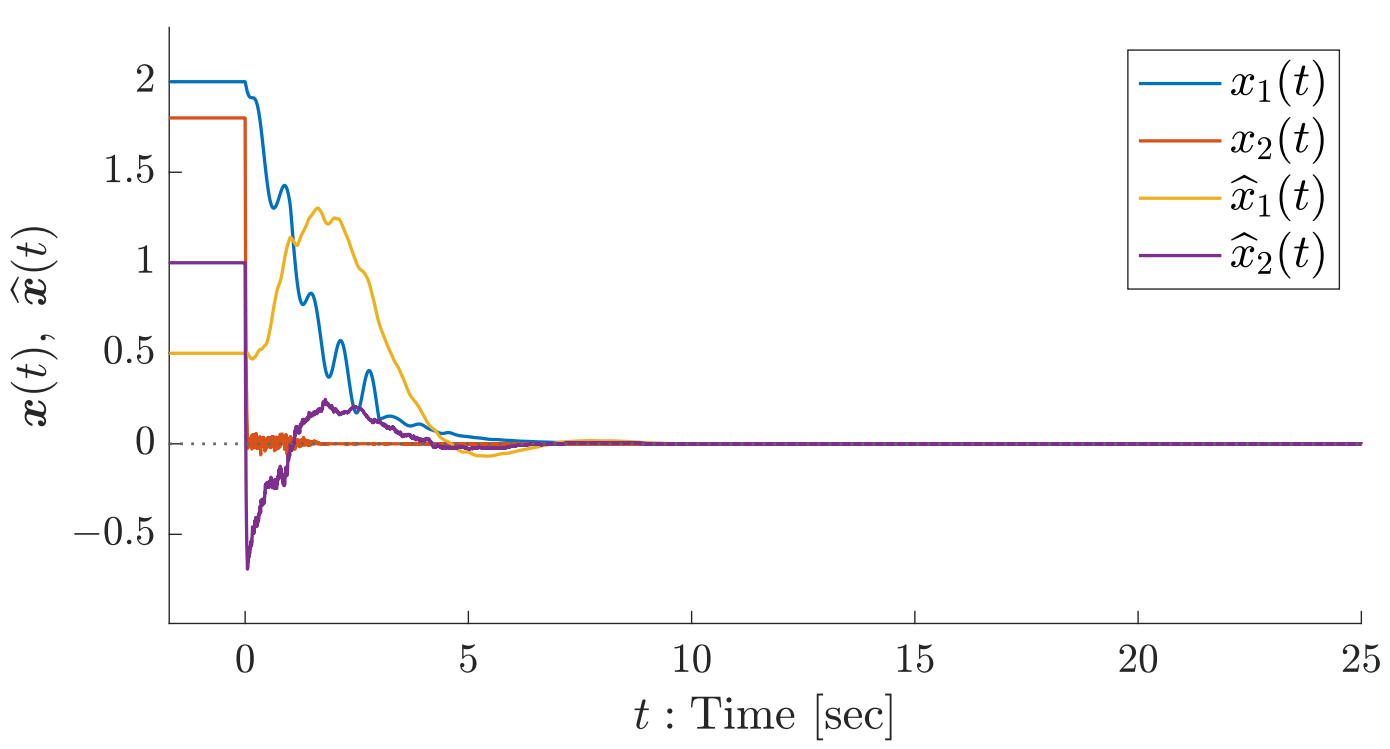

Fig. 5**:** Plot of $\boldsymbol{x}(t)$ and estimator state $\widehat{\boldsymbol{x}}(t)$

### B. DSE Design of a Nonlinear TDS

This subsection aims to demonstrate the versatility of the proposed methods in accommodating a range of applications. Specifically, assume system (1) with (42) is compensated by non-smooth static output feedback controller

$$\mathfrak{f}_1(\mathbf{u}_t(\cdot),\mathbf{y}_t(\cdot)) = -[0 \quad 30\boldsymbol{y}(t)+30|\boldsymbol{y}(t))|^c \operatorname{sign}(\boldsymbol{y}(t))]^\top \quad (45)$$

with $\mathfrak{f}_2(\mathbf{u}_t(\cdot),\mathbf{y}_t(\cdot)) = \mathfrak{f}_3(\mathbf{u}_t(\cdot)) = 0$ and the measurement output $\boldsymbol{y}(t) = \mathcal{C}\boldsymbol{x}(t) = x_2(t) \in \mathbb{R}$ as in the previous subsection, where $\mathcal{C}$ is given by (42) and $0 < c < 1$ is known.

As explained in Remark 4, equation (10) is independent of the control input/output feedbacks $\mathfrak{f}_i(\mathbf{u}_t(\cdot),\mathbf{y}_t(\cdot))$ in (1). Consequently, the error-dynamics in this subsection is **identical** to the one in the previous subsection. This implies that the estimator gains corresponding to the results in Tables I–IV are also applicable to the nonlinear system being examined.

Employing the same estimator gains corresponding to $\min \gamma = 0.5024, \mathrm{NoIs} = 15$ in Table II with the same initial conditions, disturbances and signal $\varsigma(t)$ as in the previous subsection, the simulation results are provided in Figures 5-6. Since the error-dynamics equation and initial conditions remain unchanged, the plots of $\boldsymbol{e}(t)$ and $\boldsymbol{\zeta}(t)$ are **identical** to those in Figures 2 and 4. Indeed, for the same error-dynamics, the resulting estimator gains can work for different TDSs provided that the open loop system can be denoted in the form of (1). It is evident from Figures 5–6 that our estimator achieves state and output estimation with stochastic noise $\varsigma(t)$.

### V. Conclusion and Remarks

An effective solution to the DSE design problem of system (1) with general delay structures has been developed using the proposed control concept called KDS in conjunction with the KF approach. The key idea is to decompose matrix-valued functions with $\mathcal{L}^2$ kernels using finite yet limitless numbers of compartmentalized basis functions, where a portion of the functions are approximated by the rest. This allows the basis functions to be addressed differently while enabling the construction of KFs with general $\boldsymbol{f}_i(\cdot) \in \mathcal{H}^1(\mathcal{I}_i\,; \mathbb{R}^{d_i})$ integral kernels. The KSD has been shown to work seamlessly with estimator (12) and general KF (48a) to form synthesis conditions expressed in SDP conditions with finite dimensions, where all the relevant infinite-dimensional functions are embedded in $\boldsymbol{\vartheta}(t)$ of the quadratic forms in step (54).

An interesting future study is to investigate if a result of hierarchical-like feasibility improvement can be established for our synthesis conditions with $\boldsymbol{f}_i(\cdot) \in \mathcal{H}^1(\mathcal{I}_i\,; \mathbb{R}^{d_i})$, which has been suggested by [19], [63] with orthogonal functions/polynomials. This also raises the question whether the induced conservatism of our KF approach can be reduced indefinitely by adding new functions to $\boldsymbol{f}_i(\cdot)$. Another important direction to consider is to develop a quantitative framework that can find the optimal $\boldsymbol{g}_i(\cdot)$ that may ensure reasonable feasibilities with minimal numerical burden. Finally, the proposed framework can be utilized as a blueprint for developing solutions to several open-problems, such as dissipative filtering

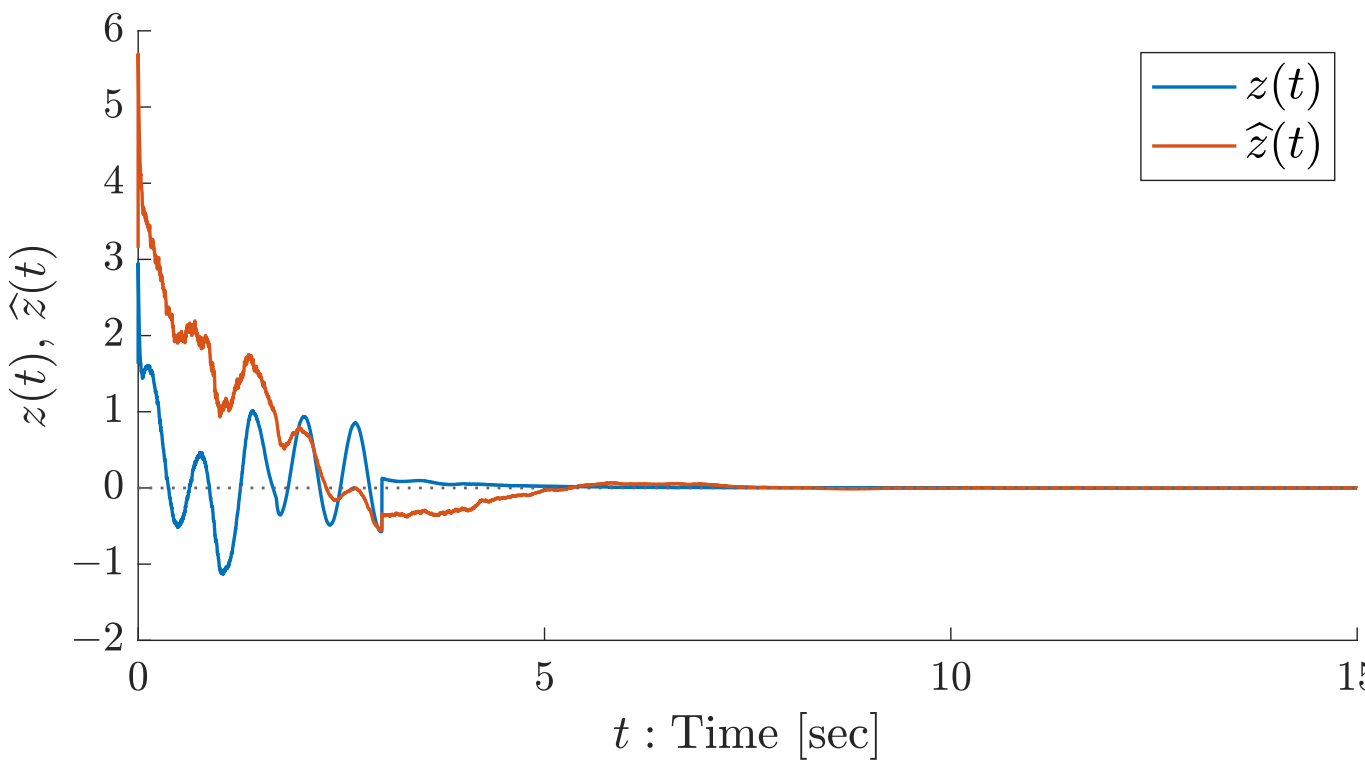

Fig. 6: Plots of regulated output $\boldsymbol{z}(t), \widehat{\boldsymbol{z}}(t)$

and dynamic output feedback of LTDSs with general delays.

## APPENDIX

We define weighted Lebesgue function space

$$\mathcal{L}^2_{\varpi}\left(\mathcal{K};\mathbb{R}^d\right) = \left\{\boldsymbol{h}(\cdot) \in \mathcal{M}\left(\mathcal{K};\mathbb{R}^d\right) : \|\boldsymbol{h}(\cdot)\|_{\varpi} < +\infty\right\} \quad (46)$$

with $d \in \mathbb{N}$ and $\|\boldsymbol{h}(\cdot)\|^2_{\varpi} = \int_{\mathcal{K}} \varpi(\tau)\boldsymbol{h}^\top(\tau)\boldsymbol{h}(\tau)\mathsf{d}\tau$, where $\varpi(\cdot) \in \mathcal{M}(\mathcal{K};\mathbb{R}_{\geq 0})$, and $\varpi(\tau)$ has countably infinite or finite number of zero values. Finally, $\mathcal{K} \subseteq \mathbb{R} \cup \{\pm\infty\}$ is Lebesgue measurable and its measure is nonzero.

*Lemma 4:* Given a sequence of sets $\mathcal{K}_i \subseteq \mathcal{K}$ and $\varpi(\cdot)$ in (46) and $\mathbb{S}^n \ni X_i \succeq 0$ with $i \in \mathbb{N}_\nu$ and $\nu \in \mathbb{N}$. Consider $\mathsf{f}_i(\cdot) \in \mathcal{L}^2_{\varpi}\left(\mathcal{K}_i;\mathbb{R}^{l_i}\right)$ and $\mathsf{g}_i(\cdot) \in \mathcal{L}^2_{\varpi}\left(\mathcal{K}_i;\mathbb{R}^{\lambda_i}\right)$ such that

$$\forall i \in \mathbb{N}_\nu, \int_{\mathcal{K}_i} \varpi(\tau) \begin{bmatrix} \mathsf{g}_i(\tau) \\ \mathsf{f}_i(\tau) \end{bmatrix} \begin{bmatrix} \mathsf{g}_i^\top(\tau) & \mathsf{f}_i^\top(\tau) \end{bmatrix} \mathsf{d}\tau \succ 0$$

with $l_i \in \mathbb{N}$ and $\lambda_i \in \mathbb{N}_0$ for all $i \in \mathbb{N}_\nu$, from which we infer $\forall i \in \mathbb{N}_\nu,\ \mathcal{F}_i = \int_{\mathcal{K}} \varpi(\tau)\mathsf{f}_i(\tau)\mathsf{f}_i^\top(\tau)\mathsf{d}\tau \succ 0$. Then

$$\begin{aligned}
&\forall \boldsymbol{x}(\cdot) \in \mathcal{L}^2_{\varpi}(\mathcal{K};\mathbb{R}^n),\ \sum_{i=1}^{\nu} \int_{\mathcal{K}_i} \varpi(\tau)\boldsymbol{x}^\top(\tau)X_i\boldsymbol{x}(\tau)\mathsf{d}\tau \\
&\geq [*]\left[\operatorname*{diag}_{i=1}^{\nu}\left(\mathcal{F}_i^{-1} \otimes X_i\right)\right]\left[\int_{\mathcal{K}_i} \varpi(\tau)\widetilde{\mathsf{F}}_i(\tau)\boldsymbol{x}(\tau)\mathsf{d}\tau\right]_{i=1}^{\nu} \\
&\quad + [*]\left[\operatorname*{diag}_{i=1}^{\nu}\left(\mathcal{E}_i^{-1} \otimes X_i\right)\right]\left[\int_{\mathcal{K}_i} \varpi(\tau)\widetilde{\mathsf{E}}_i(\tau)\boldsymbol{x}(\tau)\mathsf{d}\tau\right]_{i=1}^{\nu} \\
&\geq [*]\left[\operatorname*{diag}_{i=1}^{\nu}\left(\mathcal{F}_i^{-1} \otimes X_i\right)\right]\left[\int_{\mathcal{K}_i} \varpi(\tau)\widetilde{\mathsf{F}}_i(\tau)\boldsymbol{x}(\tau)\mathsf{d}\tau\right]_{i=1}^{\nu}, \quad (47)
\end{aligned}$$

where $\widetilde{\mathsf{F}}_i(\tau) = \mathsf{f}_i(\tau) \otimes I_n$ and $\widetilde{\mathsf{E}}_i(\tau) = \boldsymbol{\varepsilon}_i(\tau) \otimes I_n$. The rest of the symbols are defined as $\boldsymbol{\varepsilon}_i(\tau) = \mathsf{g}_i(\tau) - \mathsf{A}_i\mathsf{f}_i(\tau) \in \mathbb{R}^{\lambda_i}$ and $\mathbb{R}^{\lambda_i \times l_i} \ni \mathsf{A}_i = \int_{\mathcal{K}_i} \varpi(\tau)\mathsf{g}_i(\tau)\mathsf{f}_i^\top(\tau)\mathsf{d}\tau\mathcal{F}_i^{-1}$ and $\mathcal{E}_i = \int_{\mathcal{K}_i} \varpi(\tau)\boldsymbol{\varepsilon}_i(\tau)\boldsymbol{\varepsilon}_i^\top(\tau)\mathsf{d}\tau \succ 0$.

*Proof:* By using the inequality in [19, Eq. (17)] $\nu$ times with Lemma 1, then the first inequality in (47) is obtained. The definition of $\mathcal{F}_i$ here is the inverse matrix of the one in [19]. Furthermore, $\mathcal{E}_i^{-1} \succ 0$ in (47) is ensured by [19, Eq. (18)], from which we infer the second inequality in (47). ■

### A. Proof of Theorem 1

*Proof:* The theorem is established based on the construction of a complete type Krasovskiĭ functional

$$\begin{aligned}
\mathsf{v}(\boldsymbol{e}_t(\cdot)) &= \boldsymbol{\eta}^\top(t)\begin{bmatrix} P_1 & P_2 \\ * & P_3 \end{bmatrix}\boldsymbol{\eta}(t) \\
&\quad + \sum_{i=1}^{\nu} \int_{\mathcal{I}_i} \boldsymbol{e}^\top(t+\tau)\left[Q_i + (\tau + r_i)R_i\right]\boldsymbol{e}(t+\tau)\mathsf{d}\tau, \quad (48a)
\end{aligned}$$

where $\boldsymbol{e}_t(\cdot)$ follows the definition in (28), and

$$\begin{aligned}
&\boldsymbol{\eta}(t) = \begin{bmatrix} \boldsymbol{e}^\top(t) & \boldsymbol{\mathfrak{w}}^\top(t) \end{bmatrix}^\top, \quad F_i(\tau) = \boldsymbol{f}_i(\tau) \otimes I_n \\
&\boldsymbol{\mathfrak{w}}(t) = \left[\int_{\mathcal{I}_i}\left(\sqrt{\mathfrak{F}_i^{-1}} \otimes I_n\right)F_i(\tau)\boldsymbol{e}(t+\tau)\mathsf{d}\tau\right]_{i=1}^{\nu} \quad (48b)
\end{aligned}$$

with $\boldsymbol{f}_i(\tau)$ in (4)–(6). The definitions of $P_1, P_2, P_3, Q_i, R_i, \mathfrak{F}_i$ are given by the statements of Theorem 1. Matrices $\sqrt{\mathfrak{F}_i^{-1/2}}$ are well defined by (6) such that $\mathfrak{F}_i = \int_{\mathcal{I}_i} \boldsymbol{f}_i(\tau)\boldsymbol{f}_i^\top(\tau)\mathsf{d}\tau \succ 0$.

We first establish some crucial properties concerning the integrals in (48a). Given $\boldsymbol{\chi}(t,\tau)$ in (15), we find that

$$\begin{aligned}
&\frac{\mathsf{d}}{\mathsf{d}t}\sum_{i=1}^{\nu}\int_{\mathcal{I}_i} \boldsymbol{e}^\top(t+\tau)\ \left[Q_i + (\tau + r_i)R_i\right]\boldsymbol{e}(t+\tau)\mathsf{d}\tau \\
&= \sum_{i=1}^{\nu}[*]\ \left(Q_i + \acute{r}_i R_i\right)\boldsymbol{e}(t - r_{i-1}) - \sum_{i=1}^{\nu}[*]\,Q_i\boldsymbol{e}(t - r_i) \\
&\quad - \sum_{i=1}^{\nu}\int_{\mathcal{I}_i} \boldsymbol{e}^\top(t+\tau)R_i\boldsymbol{e}(t+\tau)\mathsf{d}\tau \\
&= \boldsymbol{\chi}^\top(t,0)\left(\mathbf{Q} + \mathbf{R}\Lambda\right)\boldsymbol{\chi}(t,0) - \boldsymbol{\chi}^\top(t,-1)\,\mathbf{Q}\,\boldsymbol{\chi}(t,-1) \\
&\quad - \sum_{i=1}^{\nu}\int_{\mathcal{I}_i} \boldsymbol{e}^\top(t+\tau)R_i\boldsymbol{e}(t+\tau)\mathsf{d}\tau, \quad (49a)
\end{aligned}$$

$$\text{and}\quad \begin{bmatrix} \boldsymbol{e}(t) \\ \boldsymbol{\chi}(t,-1) \end{bmatrix} = \begin{bmatrix} \boldsymbol{\chi}(t,0) \\ \boldsymbol{e}(t - r_\nu) \end{bmatrix} = \left[\boldsymbol{e}(t - r_i)\right]_{i=0}^{\nu} \quad (49b)$$

by the Leibniz integral rule [64] for weak derivatives, where $\mathbf{Q} = \mathsf{diag}_{i=1}^{\nu} Q_i$ and $\mathbf{R} = \mathsf{diag}_{i=1}^{\nu} R_i$ with $\Lambda$ as in (32f). Examining the $\boldsymbol{\xi}(t)$ in (24) with Lemma 1 and (4) concludes

$$\begin{aligned}
&\left[\int_{\mathcal{I}_i} F_i(\tau)\boldsymbol{e}(t+\tau)\mathsf{d}\tau\right]_{i=1}^{\nu} = \left[\int_{\mathcal{I}_i}\left[\boldsymbol{f}_i(\tau) \otimes I_n\right]\boldsymbol{e}(t+\tau)\mathsf{d}\tau\right]_{i=1}^{\nu} \\
&= \left[\int_{\mathcal{I}_i}\left(\left[\widetilde{I}_i\boldsymbol{h}_i(\tau)\right] \otimes\ I_n\right)\boldsymbol{e}(t+\tau)\mathsf{d}\tau\right]_{i=1}^{\nu} \\
&= \left[\int_{\mathcal{I}_i}\left(\widetilde{I}_i \otimes I_n\right)H_i(\tau)\boldsymbol{e}(t+\tau)\mathsf{d}\tau\right]_{i=1}^{\nu} = \left[\operatorname*{diag}_{i=1}^{\nu}\widetilde{I}_i \otimes I_n\right] \\
&\qquad \times \left[\int_{\mathcal{I}_i}\left(\sqrt{\mathfrak{H}_i}\sqrt{\mathfrak{H}_i^{-1}} \otimes I_n\right)H_i(\tau)\boldsymbol{e}(t+\tau)\mathsf{d}\tau\right]_{i=1}^{\nu} \\
&= \left[\operatorname*{diag}_{i=1}^{\nu}\widetilde{I}_i\sqrt{\mathfrak{H}_i} \otimes I_n\right]\boldsymbol{\xi}(t) \quad (50)
\end{aligned}$$

with $H_i(\tau)$ in (19), as $\boldsymbol{f}_i(\tau) = \widetilde{I}_i\boldsymbol{h}_i(\tau)$ by (4) with $\widetilde{I}_i$ in (32e). Note that $\sqrt{\mathfrak{H}_i^{-1}}$ is well defined with $\mathfrak{H}_i \succ 0$ in (8). Consider $\boldsymbol{\mathfrak{w}}(t)$ in (48b) with (50) and invoke (14a), it follows that

$$\boldsymbol{\mathfrak{w}}(t) = \left(\operatorname*{diag}_{i=1}^{\nu}\sqrt{\mathfrak{F}_i^{-1}}\widetilde{I}_i\sqrt{\mathfrak{H}_i} \otimes I_n\right)\boldsymbol{\xi}(t) = \widehat{I}\boldsymbol{\xi}(t), \quad (51)$$

where $\widehat{I}, \mathfrak{F}_i$ are given by the statements of Theorem 1. To compute the weak derivative of $\mathfrak{w}(t)$, we can first establish

$$\begin{aligned}
&\frac{\mathsf{d}}{\mathsf{d}t}\left[\int_{\mathcal{I}_i} F_i(\tau)\boldsymbol{e}(t+\tau)\mathsf{d}\tau\right]_{i=1}^{\nu} = \left[\frac{\mathsf{d}}{\mathsf{d}t}\int_{\mathcal{I}_i} F_i(\tau)\boldsymbol{e}(t+\tau)\mathsf{d}\tau\right]_{i=1}^{\nu} \\
&= \left[\int_{\mathcal{I}_i} F_i(\tau)\frac{\partial \boldsymbol{e}(t+\tau)}{\partial t}\mathsf{d}\tau\right]_{i=1}^{\nu} = \left[\int_{\mathcal{I}_i} F_i(\tau)\mathsf{d}\boldsymbol{e}(t+\tau)\right]_{i=1}^{\nu} \\
&= \Bigg[F_i(-r_{i-1})\boldsymbol{e}(t-r_{i-1}) - F_i(-r_i)\boldsymbol{e}(t-r_i) \\
&\qquad\qquad - \int_{\mathcal{I}_i}(M_i\otimes I_n)H_i(\tau)\boldsymbol{e}(t+\tau)\mathsf{d}\tau\Bigg]_{i=1}^{\nu} \\
&= \left(\mathop{\mathsf{diag}}\limits_{i=1}^{\nu} F_i(-r_{i-1})\right)\boldsymbol{\chi}(t,0) - \left(\mathop{\mathsf{diag}}\limits_{i=1}^{\nu} F_i(-r_i)\right)\boldsymbol{\chi}(t,-1) \\
&\qquad\qquad - \left(\mathop{\mathsf{diag}}\limits_{i=1}^{\nu} M_i\sqrt{\mathfrak{H}_i}\otimes I_n\right)\boldsymbol{\xi}(t)
\end{aligned}\tag{52}$$

by invoking Lemma 1 and the Leibniz integral rule [64] with the relation in (5) and $\frac{\partial \boldsymbol{e}(t+\tau)}{\partial t} = \frac{\partial \boldsymbol{e}(t+\tau)}{\partial \tau}$. In light of the identities in (51)–(52) and Lemma 1, we deduce that

$$\begin{aligned}
&\widetilde{\forall} t\geq t_0,\ \frac{\mathsf{d}\mathfrak{w}(t)}{\mathsf{d}t} = \left(\mathop{\mathsf{diag}}\limits_{i=1}^{\nu}\sqrt{\mathfrak{F}_i^{-1}}\boldsymbol{f}_i(-r_{i-1})\otimes I_n\right)\boldsymbol{\chi}(t,0) \\
&\qquad - \left(\mathop{\mathsf{diag}}\limits_{i=1}^{\nu}\sqrt{\mathfrak{F}_i^{-1}}\boldsymbol{f}_i(-r_i)\otimes I_n\right)\boldsymbol{\chi}(t,-1) \\
&\qquad - \left(\mathop{\mathsf{diag}}\limits_{i=1}^{\nu}\sqrt{\mathfrak{F}_i^{-1}}M_i\sqrt{\mathfrak{H}_i}\otimes I_n\right)\boldsymbol{\xi}(t) \\
&= \left[\mathop{\mathsf{diag}}\limits_{i=1}^{\nu}\sqrt{\mathfrak{F}_i^{-1}}\boldsymbol{f}_i(-r_{i-1})\otimes I_n \quad \mathsf{O}_{dn,n}\quad \mathsf{O}_{dn,\varkappa n}\right]\boldsymbol{\omega}(t) \\
&- \left[\mathsf{O}_{dn,n}\quad \mathop{\mathsf{diag}}\limits_{i=1}^{\nu}\sqrt{\mathfrak{F}_i^{-1}}\boldsymbol{f}_i(-r_i)\otimes I_n \quad \mathop{\mathsf{diag}}\limits_{i=1}^{\nu}\sqrt{\mathfrak{F}_i^{-1}}M_i\sqrt{\mathfrak{H}_i}\otimes I_n\right]\boldsymbol{\omega}(t) \\
&= (\mathbf{M}\otimes I_n)\boldsymbol{\omega}(t) = \left[\mathbf{M}\otimes I_n\quad \mathsf{O}_{dn,(\mu n+q)}\right]\boldsymbol{\vartheta}(t)
\end{aligned}\tag{53}$$

with $\boldsymbol{\omega}(t)$ in (26b) and $\mathbf{M}$ in (32g), where $d=\sum_{i=1}^{\nu} d_i$ with $M_i, d_i$ in (5), as $\left(\sqrt{\mathfrak{F}_i^{-1}}\otimes I_n\right)F_i(\tau) = \sqrt{\mathfrak{F}_i^{-1}}\boldsymbol{f}_i(\tau)\otimes I_n$.

With (48b)–(53) at hand, we can now differentiate (weak derivative) $\mathsf{v}(\boldsymbol{e}_t(\cdot))$ in (48a) with respect to $t$ along the trajectory of (25) with $\mathsf{s}(\boldsymbol{\zeta}(t),\boldsymbol{w}(t))$ in (30). It follows that

$$\begin{aligned}
&\widetilde{\forall} t\geq t_0,\ \dot{\mathsf{v}}(\boldsymbol{e}_t(\cdot)) - \mathsf{s}(\boldsymbol{\zeta}(t),\boldsymbol{w}(t)) = \\
&\mathsf{Sy}\left(\boldsymbol{\eta}^\top(t)\begin{bmatrix}P_1 & P_2\\ * & P_3\end{bmatrix}\dot{\boldsymbol{\eta}}(t)\right) - \begin{bmatrix}\boldsymbol{\zeta}(t)\\ \boldsymbol{w}(t)\end{bmatrix}^\top\begin{bmatrix}\widetilde{J}^\top J_1^{-1}\widetilde{J} & J_2\\ * & J_3\end{bmatrix}\begin{bmatrix}\boldsymbol{\zeta}(t)\\ \boldsymbol{w}(t)\end{bmatrix} \\
&= \boldsymbol{\vartheta}^\top(t)\,\mathsf{Sy}\Bigg(S^\top\begin{bmatrix}P_1&P_2\\ *&P_3\end{bmatrix}\begin{bmatrix}\boldsymbol{\Omega}\\ \left[\mathbf{M}\otimes I_n\quad \mathsf{O}_{dn,(\mu n+q)}\right]\end{bmatrix} \\
&\quad - \begin{bmatrix}\mathsf{O}_{(\beta n),m}\\ J_2^\top\end{bmatrix}\boldsymbol{\Sigma}\Bigg)\boldsymbol{\vartheta}(t) + \boldsymbol{\chi}^\top(t,0)\left(\mathbf{Q}+\mathbf{R}\Lambda\right)\boldsymbol{\chi}(t,0) \\
&\quad - \boldsymbol{\chi}^\top(t,-1)\,\mathbf{Q}\boldsymbol{\chi}(t,-1) - \sum_{i=1}^{\nu}\int_{\mathcal{I}_i}\boldsymbol{e}^\top(t+\tau)R_i\boldsymbol{e}(t+\tau)\mathsf{d}\tau \\
&\quad - \boldsymbol{w}^\top(t)J_3\boldsymbol{w}(t) - \boldsymbol{\vartheta}^\top(t)\boldsymbol{\Sigma}^\top\widetilde{J}^\top J_1^{-1}\widetilde{J}\boldsymbol{\Sigma}\boldsymbol{\vartheta}(t),
\end{aligned}\tag{54}$$

since $\boldsymbol{\eta}(t)=S\boldsymbol{\vartheta}(t)$ with $\boldsymbol{\eta}(t)$ in (48b) and $\boldsymbol{\vartheta}(t)$ in (26a) and $S$ in (32c), where $\boldsymbol{\Omega},\boldsymbol{\Sigma}$ are given by Theorem 1. We note that the structure of $S$ follows from $\mathfrak{w}(t) = \widehat{I}\boldsymbol{\xi}(t)$ in (51) once we realize that vectors $\boldsymbol{e}(t)$, $\boldsymbol{\xi}(t)$ are part of $\boldsymbol{\vartheta}(t)$ in (26a) and vectors $\boldsymbol{e}(t)$, $\mathfrak{w}(t)$ constitute $\boldsymbol{\eta}(t)$ in (48b).

Assume that (31b) is true. Apply (47) with the parameters in (8)–(10) such that $\varpi(\tau)=1, \mathsf{g}_i(\tau)=\boldsymbol{\phi}_i(\tau)$, $\mathsf{f}_i(\tau)=\boldsymbol{h}_i(\tau)$, $\mathsf{\varepsilon}_i(\tau)=\boldsymbol{\varepsilon}_i(\tau)$ and $X_i=R_i$, $\mathcal{F}_i=\mathfrak{H}_i$ and $\mathcal{E}_i=\mathfrak{E}_i$ and $\boldsymbol{x}(\tau)\leftarrow\boldsymbol{e}(t+\tau)$ to the quadratic integral sum in (54), then

$$\begin{aligned}
&\sum_{i=1}^{\nu}\int_{\mathcal{I}_i}\boldsymbol{e}^\top(t+\tau)R_i\boldsymbol{e}(t+\tau)\mathsf{d}\tau \geq \\
&\boldsymbol{\xi}^\top(t)\left[\mathop{\mathsf{diag}}\limits_{i=1}^{\nu}\left(I_{\varkappa_i}\otimes R_i\right)\right]\boldsymbol{\xi}(t) + \boldsymbol{\mathfrak{a}}^\top(t)\left[\mathop{\mathsf{diag}}\limits_{i=1}^{\nu}\left(I_{\mu_i}\otimes R_i\right)\right]\boldsymbol{\mathfrak{a}}(t)
\end{aligned}\tag{55}$$

with $\boldsymbol{\xi}(t)$ and $\boldsymbol{\mathfrak{a}}(t)$ as in (24) and $\varkappa_i,\mu_i$ in Proposition 1. Note that $\mathfrak{H}_i^{-1}\succ 0$ and $\mathfrak{E}_i^{-1}\succ 0$ in (55) were factorized as the products of their unique square roots that are placed in $\boldsymbol{\xi}(t),\boldsymbol{\mathfrak{a}}(t)$, respectively, owning to $\mathfrak{H}_i^{-1}\otimes R_i = [*]\left(I_{\varkappa_i}\otimes R_i\right)\left(\mathfrak{H}_i^{-1/2}\otimes I_n\right)$ and $\mathfrak{E}_i^{-1}\otimes R_i = [*]\left(I_{\mu_i}\otimes R_i\right)\left(\mathfrak{E}_i^{-1/2}\otimes I_n\right)$ by (14a). Next, utilize (55) on (54), it follows from the properties of quadratic forms with (49b) and $\boldsymbol{\vartheta}(t)$ in (26a) that

$$\begin{aligned}
&\widetilde{\forall} t\geq t_0,\ \dot{\mathsf{v}}(\boldsymbol{e}_t(\cdot)) - \mathsf{s}(\boldsymbol{\zeta}(t),\boldsymbol{w}(t)) \\
&\qquad\qquad \leq \boldsymbol{\vartheta}^\top(t)\left(\boldsymbol{\Psi} - \boldsymbol{\Sigma}^\top\widetilde{J}^\top J_1^{-1}\widetilde{J}\boldsymbol{\Sigma}\right)\boldsymbol{\vartheta}(t)
\end{aligned}$$

with $\boldsymbol{\Psi}$ in (32b). Appealing to the properties of quadratic form again, we can conclude that there exists $\epsilon_3>0$ such that

$$\widetilde{\forall} t\geq t_0,\ \dot{\mathsf{v}}(\boldsymbol{e}_t(\cdot)) - \mathsf{s}(\boldsymbol{\zeta}(t),\boldsymbol{w}(t)) \leq -\epsilon_3\left\|\boldsymbol{e}(t)\right\|^2 \tag{56}$$

for any $\boldsymbol{\psi}(\cdot)\in\mathcal{C}(\mathcal{J};\mathbb{R}^n)$ in (25), provided that (31b) and

$$\boldsymbol{\Psi}-\boldsymbol{\Sigma}^\top\widetilde{J}^\top J_1^{-1}\widetilde{J}\boldsymbol{\Sigma}\prec 0 \tag{57}$$

hold. Given the positions of $J_2, J_3, \widetilde{J}$ in $\boldsymbol{\Psi}$ and $J_1$ in (57) and $\boldsymbol{w}(t)$ in $\boldsymbol{\vartheta}(t)$, it follows from (56)–(57) that

$$\exists\epsilon_3>0,\ \widetilde{\forall}t\geq t_0,\ \dot{\mathsf{v}}(\boldsymbol{e}_t(\cdot))\leq-\epsilon_3\left\|\boldsymbol{e}(t)\right\|^2 \tag{58}$$

for any $\boldsymbol{\psi}(\cdot)\in\mathcal{C}(\mathcal{J};\mathbb{R}^n)$ in (25) with $\boldsymbol{w}(t)\equiv\mathbf{0}_q$. Note that $\boldsymbol{e}_t(\cdot)$ in (58) follows the definition in (27b). Hence $\mathsf{v}(\boldsymbol{e}_t(\cdot))$ in (48a) satisfies (27b)–(28) if (31b) and (57) are feasible. Finally, applying the Schur complement to (57) with (31b) and $J_1^{-1}\prec 0$ yields (31c). Thus we have shown that $\mathsf{v}(\boldsymbol{e}_t(\cdot))$ satisfies (27b)–(28) with some $\epsilon_3>0$ if (31b)–(31c) are feasible.

Now we turn to prove that $\mathsf{v}(\boldsymbol{e}_t(\cdot))$ in (48a) satisfies (27a) for some $\epsilon_1;\epsilon_2>0$ if (31a)–(31b) are feasible. Consider $\mathsf{v}(\boldsymbol{e}_t(\cdot))$ at $t=t_0$ with $\boldsymbol{\eta}(t)$ as in (48b), we find that

$$\begin{aligned}
&\exists\lambda>0,\ \forall\boldsymbol{\psi}(\cdot)\in\mathcal{C}\left(\mathcal{J};\mathbb{R}^n\right),\ \mathsf{v}(\boldsymbol{e}_{t_0}(\cdot)) = \mathsf{v}(\boldsymbol{\psi}(\cdot)) \\
&\leq \boldsymbol{\eta}^\top(t_0)\lambda\boldsymbol{\eta}(t_0) + \int_{-r}^0\boldsymbol{\psi}^\top(\tau)\lambda\boldsymbol{\psi}(\tau)\mathsf{d}\tau \\
&\leq \lambda\left\|\boldsymbol{\psi}(0)\right\|^2 + \lambda r\|\boldsymbol{\psi}(\cdot)\|_\infty^2 + [*]\left[\mathop{\mathsf{diag}}\limits_{i=1}^{\nu}\left(I_{d_i}\otimes\lambda I_n\right)\right]\mathfrak{w}(t_0) \\
&\leq \lambda\left\|\boldsymbol{\psi}(0)\right\|^2+\lambda r\left\|\boldsymbol{\psi}(\cdot)\right\|_\infty^2 + \sum_{i=1}^{\nu}\int_{\mathcal{I}_i}\boldsymbol{\psi}^\top(\tau)\lambda\boldsymbol{\psi}(\tau)\mathsf{d}\tau \\
&\leq(\lambda+\lambda r)\left\|\boldsymbol{\psi}(\cdot)\right\|_\infty^2+\lambda\int_{-r}^0\boldsymbol{\psi}^\top(\tau)\boldsymbol{\psi}(\tau)\mathsf{d}\tau \\
&\leq(\lambda+2\lambda r)\left\|\boldsymbol{\psi}(\cdot)\right\|_\infty^2,
\end{aligned}\tag{59}$$

following from the property $\forall X\in\mathbb{S}^n$, $\exists\lambda>0$, $\forall\mathbf{x}\in\mathbb{R}^n\setminus\{\mathbf{0}_n\}, \mathbf{x}^\top\left(\lambda I_n - X\right)\mathbf{x}>0$ and the second inequality in (47) with $\varpi(\tau)=1$, $\mathsf{f}_i(\tau)=\sqrt{\mathfrak{F}_i^{-1}}\boldsymbol{f}_i(\tau)$ and $X_i=\lambda I_n$. Now applying the second inequality in (47) again to (48a) with $\varpi(\tau)=1$, $\mathcal{K}_i=\mathcal{I}_i$, $\mathsf{f}_i(\tau)=\sqrt{\mathfrak{F}_i^{-1}}\boldsymbol{f}_i(\tau)$ and $X_i=Q_i$ yields

$$\sum_{i=1}^{\nu}\int_{\mathcal{I}_i}[*]Q_i\boldsymbol{e}(t+\tau)\mathsf{d}\tau \geq [*]\left[\mathop{\mathsf{diag}}_{i=1}^{\nu}\left(I_{d_i}\otimes Q_i\right)\right]\mathfrak{w}(t) \quad (60)$$

with $\mathfrak{w}(t)$ in (48b), provided that (31b) holds. Having established the relations in (60), (31b) and (59), it indicates that (48a) satisfies (27a) for some $\epsilon_1;\epsilon_2>0$ if (31a)–(31b) hold.

In conclusion, we have shown that feasible solutions to (31a)–(31c) imply the existence of the KF in (48a) and some constants $\epsilon_1;\epsilon_2;\epsilon_3>0$ satisfying the stability criteria (27a)–(27b) and the dissipativity condition (28). ■

*Remark 10:* All the DD kernels in (3) must be present in $\boldsymbol{g}_i(\cdot)$ by (4) in Proposition 1. Meanwhile, an unlimited number of WDLIFs can be added to $\boldsymbol{f}_i(\cdot)\in\mathscr{H}^1(\mathcal{I}_i;\mathbb{R}^{d_i})$ in (4)–(6) and (48b) to increase the generality of KF (48a). This in turn improves the feasibility of (31a)–(31c), even if no functions inside of $\boldsymbol{f}_i(\cdot)$ are included by any DD integral matrix in (3). The numerical complexity of Theorem 1 is $\mathcal{O}(d^2n^2)$, depending on $d_i=\mathsf{dim}(\boldsymbol{f}_i(\tau))$. Thus, the choice of $\boldsymbol{g}_i(\cdot)$ should balance both feasibility and computational complexity. Finally, the equations in (5) are specifically formulated to make sure that both $\boldsymbol{\Omega}$ and $\mathbf{M}$ in (54) are constructible within a quadratic form via $\boldsymbol{\vartheta}(t)$ in (26a), given the structure of $\mathfrak{w}(t)$ in (48b) and the expressions of $\frac{\mathsf{d}}{\mathsf{d}t}\mathfrak{w}(t)$ in (53).

### B. Conservatism Analysis

The KF in (48b) is a parameterization of complete KF

$$\begin{aligned}V(\boldsymbol{e}_t(\cdot))&=\boldsymbol{e}^\top(t)P_1\boldsymbol{e}(t)+2\sum_{i=1}^{\nu}\boldsymbol{e}^\top(t)\int_{\mathcal{I}_i}\widetilde{P}_i^2(\tau)\boldsymbol{e}(t+\tau)\mathsf{d}\tau\\&+\sum_{i=1}^{\nu}\sum_{j=1}^{\nu}\int_{\mathcal{I}_i}\int_{\mathcal{I}_j}\boldsymbol{e}^\top(t+\theta)\widetilde{P}_{i,j}^3(\tau,\theta)\boldsymbol{e}(t+\tau)\mathsf{d}\tau\mathsf{d}\theta\\&+\sum_{i=1}^{\nu}\int_{\mathcal{I}_i}\boldsymbol{e}^\top(t+\tau)\left(Q_i+(\tau+r_i)R_i\right)\boldsymbol{e}(t+\tau)\mathsf{d}\tau \quad (61)\end{aligned}$$

for LTDS $\dot{\boldsymbol{e}}(t)=\sum_{i=0}^{\nu}A_i\boldsymbol{e}(t-r_i)$ using basis functions $\sqrt{\mathfrak{F}_i^{-1}}\boldsymbol{f}_i(\tau)$, which can well balance generality and computational complexity. The structure of (61) is equivalent to the complete KF in [22] by the properties of integrals, where the matrices are defined as $P_1\in\mathbb{S}^n$, $\widetilde{P}_i^2(\cdot)\in\mathscr{H}^1(\mathcal{I}_i;\mathbb{R}^{n\times n})$ and $\widetilde{P}_{i,j}^3(\cdot,\cdot)\in\mathscr{H}^1(\mathcal{J}^2;\mathbb{R}^{n\times n})$ satisfying $\forall(\tau,\theta)\in\mathcal{J}^2$, $\widetilde{P}_{i,j}^3(\tau,\theta)=\widetilde{P}_{j,i}^{3,\top}(\theta,\tau)$. Now $\widetilde{P}_i^2(\cdot),\widetilde{P}_{i,j}^3(\cdot,\cdot)$ are infinite-dimensional unknowns, which are difficult to compute numerically. By employing subspace parameterizations via $\sqrt{\mathfrak{F}_i^{-1}}\boldsymbol{f}_i(\tau)$ with $\mathfrak{F}_i$ and $\boldsymbol{f}_i(\cdot)\in\mathscr{H}^1(\mathcal{I}_i;\mathbb{R}^{d_i})$ in Theorem 1 and Proposition 1, we can set

$$\widetilde{P}_i^2(\tau)=\widehat{P}_i^2\left[\sqrt{\mathfrak{F}_i^{-1}}\boldsymbol{f}_i(\tau)\otimes I_n\right],\quad \widehat{P}_i^2\in\mathbb{R}^{n\times d_in}$$

$$\widetilde{P}_{i,j}^3(\tau,\theta)=\left[\sqrt{\mathfrak{F}_j^{-1}}\boldsymbol{f}_j(\theta)\otimes I_n\right]^\top\widehat{P}_{i,j}^3\left[\sqrt{\mathfrak{F}_i^{-1}}\boldsymbol{f}_i(\tau)\otimes I_n\right]$$

with $\widehat{P}_{i,j}^3\in\mathbb{R}^{d_in\times d_jn}$ satisfying $\widehat{P}_{i,j}^3=\widehat{P}_{j,i}^{3,\top}$. Thus the complete KF in (61) becomes our KF in (48a) if we let $P_2=\left[\!\left[\widehat{P}_i^2\right]\!\right]_{i=1}^{\nu}$ and $P_3=\left[\!\left[\left[\widehat{P}_{i,j}^3\right]_{i=1}^{\nu}\right]\!\right]_{j=1}^{\nu}$ in (48a) since

$$\begin{aligned}&\sum_{i=1}^{\nu}\boldsymbol{e}^\top(t)\int_{\mathcal{I}_i}\widetilde{P}_i^2(\tau)\boldsymbol{e}(t+\tau)\mathsf{d}\tau\\&=\boldsymbol{e}^\top(t)\sum_{i=1}^{\nu}\widehat{P}_i^2\int_{\mathcal{I}_i}\left[\sqrt{\mathfrak{F}_i^{-1}}\boldsymbol{f}_i(\tau)\otimes I_n\right]\boldsymbol{e}(t+\tau)\mathsf{d}\tau\\&=\boldsymbol{e}^\top(t)\left[\!\left[\widehat{P}_i^2\right]\!\right]_{i=1}^{\nu}\mathfrak{w}(t)=\boldsymbol{e}^\top(t)P_2\mathfrak{w}(t)=\mathfrak{w}^\top(t)P_2^\top\boldsymbol{e}(t)\end{aligned}$$

and

$$\begin{aligned}&\sum_{i=1}^{\nu}\sum_{j=1}^{\nu}\int_{\mathcal{I}_i}\int_{\mathcal{I}_j}\boldsymbol{e}^\top(t+\theta)\widetilde{P}_{i,j}^3(\tau,\theta)\boldsymbol{e}(t+\tau)\mathsf{d}\tau\mathsf{d}\theta\\&=\sum_{i=1}^{\nu}\sum_{j=1}^{\nu}\int_{\mathcal{I}_j}\boldsymbol{e}^\top(t+\theta)\left[\sqrt{\mathfrak{F}_j^{-1}}\boldsymbol{f}_j(\theta)\otimes I_n\right]^\top\mathsf{d}\theta\,\widehat{P}_{i,j}^3\times\\&\int_{\mathcal{I}_i}\left[\sqrt{\mathfrak{F}_i^{-1}}\boldsymbol{f}_i(\tau)\otimes I_n\right]\boldsymbol{e}(t+\tau)\mathsf{d}\tau=\mathfrak{w}^\top(t)\left[\!\left[\left[\widehat{P}_{i,j}^3\right]_{i=1}^{\nu}\right]\!\right]_{j=1}^{\nu}\mathfrak{w}(t)\end{aligned}$$

with $\mathfrak{w}(t)$ in (48b), based on the definition of quadratic forms and $\left[\int_{\mathcal{I}_i}F_i(\tau)\boldsymbol{e}(t+\tau)\mathsf{d}\tau\right]_{i=1}^{\nu}=\left[\int_{\mathcal{I}_j}F_j(\theta)\boldsymbol{e}(t+\theta)\mathsf{d}\theta\right]_{j=1}^{\nu}$.

The feasibility of (31a)–(31c) is influenced by the generality of (48a) which in turn hinges on the scope of $\boldsymbol{f}_i(\cdot)\in\mathscr{H}^1\left(\mathcal{I}_i;\mathbb{R}^{d_i}\right)$. Thus, Theorem 1 is ensured to be nonconservative. In terms of stability analysis, wherein all estimator gains are zero, the KFs and stability conditions in [19], [26], [37] for LTDSs are the special cases of our KF (48a) and Theorem 1 with particular $\boldsymbol{f}_1(\cdot)$, $\nu=1$, which is far more conservative than $\boldsymbol{f}_i(\cdot)\in\mathscr{H}^1\left(\mathcal{I}_i;\mathbb{R}^{d_i}\right)$ in Proposition 1. As it has been shown that the approaches in [19], [26], [37] can detect the delay margins of numerous challenging TDS examples, the effectiveness of our method is guaranteed, as Theorem 1 cannot be more conservative than its special cases in [19], [26], [37]. Notably, recent discoveries in [63], [65] reveal that the SDP conditions constructed with Legendre polynomial $\boldsymbol{f}_1(\tau)=\boldsymbol{\ell}_d(\tau)$ are both sufficient and necessary for the existence of complete KF for $\dot{\boldsymbol{x}}(t)=A_1\boldsymbol{x}(t)+A_2\boldsymbol{x}(t-r)$, given a sufficiently large $d$. Finally, the inequality in (47), which is utilized in (55) and (60), is derived based on the least-squares approximations [39, page 182]. This ensures that the induced conservatism at steps (55) and (60) is minimized.

### C. Proof of Theorem 2

*Proof:* First, we rewrite $\mathsf{Sy}\left(\mathbf{P}^\top\boldsymbol{\Pi}\right)+\boldsymbol{\Phi}\prec 0$ in (31c) as

$$\mathsf{Sy}\left(\mathbf{P}^\top\boldsymbol{\Pi}\right)+\boldsymbol{\Phi}=[*]\begin{bmatrix}\mathsf{O}_n & \mathbf{P}\\ * & \boldsymbol{\Phi}\end{bmatrix}\begin{bmatrix}\boldsymbol{\Pi}\\ I_{\beta n+q+m}\end{bmatrix}\prec 0 \quad (62)$$

with $\beta$ in (25). We need to construct an extra inequality to apply Lemma 3 to (62), as there are two inequalities in (34b). By inspecting the bottom-right corner of $\boldsymbol{\Phi}$, we see that

$$\Upsilon^\top\begin{bmatrix}\mathsf{O}_n & \mathbf{P}\\ * & \boldsymbol{\Phi}\end{bmatrix}\Upsilon=\begin{bmatrix}-J_3-\mathsf{Sy}(D_2^\top J_2) & D_2^\top\widetilde{J}\\ * & J_1\end{bmatrix}\prec 0, \quad (63)$$

where $\Upsilon^\top:=\begin{bmatrix}\mathsf{O}_{(q+m),(n+\beta n)} & I_{q+m}\end{bmatrix}$. Obviously, the inequality in (63) is implied by (62) or (31c), as the matrix in (63) is the $2\times 2$ block matrix at the bottom-right corner of $\mathsf{Sy}\left(\mathbf{P}^\top\boldsymbol{\Pi}\right)+\boldsymbol{\Phi}$ and $\boldsymbol{\Phi}$. Furthermore, the following identities

$$\begin{aligned}&\begin{bmatrix}-I_n & \boldsymbol{\Pi}\end{bmatrix}\begin{bmatrix}-I_n & \boldsymbol{\Pi}\end{bmatrix}_\perp=\begin{bmatrix}-I_n & \boldsymbol{\Pi}\end{bmatrix}\begin{bmatrix}\boldsymbol{\Pi}\\ I_{\beta n+q+m}\end{bmatrix}=\mathsf{O}_{n,(\beta n+q+m)}\\&\widetilde{\Upsilon}^\top=\begin{bmatrix}I_{n+\beta n}\\ \mathsf{O}_{(q+m),(n+\beta n)}\end{bmatrix},\ \widetilde{\Upsilon}_\perp=\Upsilon=\begin{bmatrix}\mathsf{O}_{(n+\beta n),(q+m)}\\ I_{q+m}\end{bmatrix}, \quad (64)\\&\widetilde{\Upsilon}\Upsilon=\mathsf{O}_{(n+\beta n),(q+m)}\end{aligned}$$

indicate $\dim(\mathrm{Ker}\begin{bmatrix}-I_n & \mathbf{\Pi}\end{bmatrix}) = \mathrm{rank}\left(\begin{bmatrix}\mathbf{\Pi}^\top & I_{\beta n+q+m}\end{bmatrix}\right) = \beta n + q + m$ and $\dim(\mathrm{Ker}\,\widetilde{\Upsilon}) = \mathrm{rank}\,\Upsilon = q + m$ by the rank nullity theorem [46], from which we know that $P := \begin{bmatrix}-I_n & \mathbf{\Pi}\end{bmatrix}$ and $Q := \widetilde{\Upsilon}$ can be utilized with Lemma 3.

By applying Lemma 3 to (62)–(63) with (64), it shows (62)–(63) hold if and only if $\exists \mathbf{W} \in \mathbb{R}^{(n+\beta n)\times n}$ such that

$$\mathsf{Sy}\left(\begin{bmatrix} I_{n+\beta n} \\ \mathsf{O}_{(q+m),(n+\beta n)}\end{bmatrix} \mathbf{W} \begin{bmatrix}-I_n & \mathbf{\Pi}\end{bmatrix}\right) + \begin{bmatrix}\mathsf{O}_n & \mathbf{P} \\ * & \mathbf{\Phi}\end{bmatrix} \prec 0. \quad (65)$$

Now the inequality in (65) is still bilinear due to the product between $\mathbf{W}$ and $\mathbf{\Pi}$. To address this problem, we can set

$$\mathbf{W} = \begin{bmatrix} W \\ [\alpha_i W]_{i=1}^{\beta}\end{bmatrix}, \quad W \in \mathbb{S}^n, \quad (66)$$

with given scalars $\{\alpha_i\}_{i=1}^{\beta} \subset \mathbb{R}$. With (66), Eq. (65) becomes

$$\mathbf{\Theta} = \mathsf{Sy}\left(\begin{bmatrix} W \\ [\alpha_i W]_{i=1}^{\beta} \\ \mathsf{O}_{(q+m),n}\end{bmatrix}\begin{bmatrix}-I_n & \mathbf{\Pi}\end{bmatrix}\right) + \begin{bmatrix}\mathsf{O}_n & \mathbf{P} \\ * & \mathbf{\Phi}\end{bmatrix} \prec 0, \quad (67)$$

from which we infer (62). Note that (67) is only a sufficient condition for (62) or (31c) due to the structural constraints in (66). To convexify the bilinear products in (67), we can set

$$\begin{aligned}\widehat{\mathbf{\Pi}} = W\mathbf{\Pi} &= \begin{bmatrix} W\mathbf{A} + W\mathbf{L}_1\left[(I_\beta \otimes \mathcal{C}) \oplus \mathsf{O}_q\right] & \mathsf{O}_{n,m}\end{bmatrix} \\ &= \begin{bmatrix} W\mathbf{A} + \mathbf{U}\left[(I_\beta \otimes \mathcal{C}) \oplus \mathsf{O}_q\right] & \mathsf{O}_{n,m}\end{bmatrix} \quad (68)\end{aligned}$$

with $\mathbf{U}$ in (36a), where $U_0 = WL_0$, $U_i = WL_i$, $\widehat{U}_i = W\widehat{L}_i$. Note that $W^{-1} \succ 0$ is implied by (67), as $-W - W^\top = -2W$ is the sole element in the first diagonal block of $\mathbf{\Theta}$.

With (68) at hand, Eq. (35) follows from the inequality in (67) that implies (31c). As such, (31a)–(31b) with (35) imply (31a)–(31c), which in turn entails (27a)–(27b) and (28). ■

*Remark 11:* Although step (66) can introduce conservatism, the structure of $P_1, P_2$ in our KF (48a) remains unchanged. Such an arrangement is less conservative than directly letting $P_2 = [\![ b_i P_1 ]\!]_{i=1}^{d}$ with known $\{b_i\}_{i=1}^{d} \subset \mathbb{R}$ to convexify (31c). We also note that the use of slack variables in (65) does not enhance its feasibility compared to (31c) in Theorem 1.